# A Galerkin Isogeometric Method for Karhunen-Loève Approximation of Random Fields[☆]

Sharif Rahman[1]

*College of Engineering, The University of Iowa, Iowa City, Iowa 52242, U.S.A.*

## Abstract

This paper marks the debut of a Galerkin isogeometric method for solving a Fredholm integral eigenvalue problem, enabling random field discretization by means of the Karhunen-Loève expansion. The method involves a Galerkin projection onto a finite-dimensional subspace of a Hilbert space, basis splines (B-splines) and non-uniform rational B-splines (NURBS) spanning the subspace, and standard methods of eigensolutions. Compared with the existing Galerkin methods, such as the finite-element and mesh-free methods, the NURBS-based isogeometric method upholds exact geometrical representation of the physical or computational domain and exploits regularity of basis functions delivering globally smooth eigensolutions. Therefore, the introduction of the isogeometric method for random field discretization is not only new; it also offers a few computational advantages over existing methods. In the big picture, the use of NURBS for random field discretization enriches the isogeometric paradigm. As a result, an uncertainty quantification pipeline of the future can be envisioned where geometric modeling, stress analysis, and stochastic simulation are all integrated using the same building blocks of NURBS. Three numerical examples, including a three-dimensional random field discretization problem, illustrate the accuracy and convergence properties of the isogeometric method for obtaining eigensolutions.

*Keywords:* B-splines, NURBS, Fredholm integral eigenvalue problem, Hilbert-Schmidt operator, uncertainty quantification.

## 1. Introduction

Many uncertainty quantification problems in engineering and applied sciences require modeling spatial variability of random input parameters. For instance, the tensile and fracture toughness properties of engineering materials, the size and shape characteristics of mechanical components, and the wind and snow loads in structural systems all exhibit randomness that varies not only from sample to sample, but also from point to point in their respective domains. Therefore, random field treatment of spatial varying randomness is a vital ingredient in computational analysis. Loosely speaking, a random field represents a random quantity at each point of the domain and, therefore, engenders an infinite number of random variables. In practice, though, the number of random variables must be finite and manageable but also large enough to ensure an optimal or accurate approximation of the original random field. This process is often referred to as random field discretization.

A number of methods and approaches are available for random field discretization. For brevity, this paper will not perform a comprehensive review, but will direct readers to a paper by Betz *et al.* [1], which provides a summary of existing works, including many references cited therein. A popular approach, known by the name of Karhunen-Loève (KL) expansion [2, 3, 4], entails spectral decomposition of the covariance

---
[☆]Grant sponsor: U.S. National Science Foundation; Grant No. CMMI-1607398.
   *Email address:* sharif-rahman@uiowa.edu (Sharif Rahman)
   [1]Professor.

*Preprint submitted to Computer Methods in Applied Mechanics and Engineering*      *April 16, 2018*

function, leading to an infinite series consisting of deterministic functions of space and uncorrelated random variables. The expansion is well known with diverse applications in engineering and applied sciences [5]. However, the KL expansion mandates solution of a Fredholm integral eigenvalue problem [6], which is not an easy task in general. Analytical solutions are available only when the covariance function has simpler functional forms, such as exponential or linear functions, and/or the problem domain is rectangular. For arbitrary covariance functions or arbitrary domains in two or three dimensions, numerical methods are often needed to solve the eigenvalue problem approximately.

For numerical solution of the integral eigenvalue problem, a well-known method is the Galerkin finite-element method (FEM) employed by Ghanem and Spanos [7] in the 1990s. Roughly speaking, the finite-element solution consists of a variational formulation and function spaces defined by its basis functions [8]. These basis functions are described by local representations via finite elements, resulting in a mesh or grid, which constitutes a non-overlapping decomposition of the computational domain into elementary shapes, such as triangles or tetrahedra and quadrilaterals or hexahedra. However, for mechanical systems with complex geometry, a finite-element mesh is often created from a computer-aided design (CAD) model, where the former is an approximation of the latter. Therefore, an additional source of imprecision is embedded in the FEM-based eigensolution. Another Galerkin approach, which sidesteps the need for element-wise decomposition, is the meshless or mesh-free method, especially the element-free Galerkin method [9], upon which Rahman and Xu [10, 11] capitalize for the solution of the integral eigenvalue problem. The fundamental aspects of both FEM and the mesh-free method are identical as they are rooted in the same Galerkin formulation, but the function spaces and their basis functions are different: in FEM, the basis functions are interpolatory polynomials with $C^0$-continuity across element boundaries, whereas in the mesh-free method, the basis functions are non-interpolatory rational functions with at least $C^1$-continuity everywhere. In consequence, the approximate eigenfunctions of the KL expansion obtained by the mesh-free method are usually globally smoother than those derived from FEM. Nonetheless, as in FEM, the link between the mesh-free method and CAD geometry is, at best, tenuous [12]. Indeed, FEM or the mesh-free method may never faithfully replicate the CAD geometry. More importantly, for complex engineering designs, generating a high-quality finite-element mesh or mesh-free discretization from the CAD geometry is more formidable than performing the analysis. This is the principal motivation behind replacing finite-element- or mesh-free-generated basis functions with CAD-generated basis functions for solving the integral eigenvalue problem directly, leading to effective random field discretization.

This paper presents a Galerkin isogeometric method for solving the integral eigenvalue problem stemming from the KL expansion of a random field with an arbitrary covariance function and an arbitrary computational domain in three dimensions. The method entails performing a Galerkin discretization of the integral eigenvalue problem, formulation of the associated matrix eigenvalue problem by constructing the isogeometric function spaces spanned by basis splines (B-splines) and non-uniform rational B-splines (NURBS), and solution of the resultant matrix eigenvalue problem by standard methods. The paper is organized as follows. A brief exposition of NURBS paraphernalia and isogeometric concept is given in Section 2. Section 3 formally defines a random field and its KL expansion, followed by truncation of the KL expansion and a description of associated error measures. The limitation of the KL expansion is also discussed. Section 4 presents the proposed isogeometric method for solving the integral eigenvalue problem. The properties and construction of system matrices involved in the matrix eigenvalue problem are explained. The results from three numerical examples of increasing dimensions and hence complexity are reported in Section 5 and Appendix A. Section 6 discusses future work. Finally, conclusions are drawn in Section 7.

## 2. Isogeometric Analysis

Let $\mathbb{N} := \{1, 2, \ldots\}$, $\mathbb{N}_0 := \mathbb{N} \cup \{0\}$, $\mathbb{R} := (-\infty, +\infty)$, $\mathbb{R}_0^+ := [0, +\infty)$, and $\mathbb{R}^+ := (0, +\infty)$ represent the sets of positive integer (natural), non-negative integer, real, non-negative real, and positive real numbers, respectively. Denote by $d$ the dimension of the physical or computational domain $\mathcal{D}$ of a geometrical object, which can be a curve, surface, and solid in $\mathbb{R}^d$. In this work, $d = 1, 2, 3$, and $\mathcal{D} \subset \mathbb{R}^d$ is a bounded domain. These standard notations will be used throughout the paper.

The isogeometric analysis (IGA) employs basis functions from CAD, such as B-splines and NURBS, directly in computational analysis. Hughes *et al.* [13] were the first to propose the isogeometric paradigm and its computational framework. For the paper to be self-contained, a brief summary of NURBS-based IGA is presented here. The description is restricted to geometries modeled as a single patch. However, for NURBS-based IGA, it is sometimes necessary to represent the physical or computational domain by a multi-patch geometric model, for example, when analyzing multiply-connected domains. The multi-patch geometries were not considered in this work.

### 2.1. Knot Vectors

Consider a $d$-dimensional cartesian coordinate system in the parametric domain $\hat{\mathcal{D}} = [0,1]^d$, where an arbitrary point has coordinate $\boldsymbol{\xi} = (\xi_1, \ldots, \xi_d)$. For the coordinate direction $k$, where $k = 1, \ldots, d$, define a positive integer $n_k \in \mathbb{N}$ and a non-negative integer $p_k \in \mathbb{N}_0$, representing the total number of basis functions and polynomial degree, respectively. [2] Given $n_k$ and $p_k$, introduce on the parametric interval $[0,1] \subset \mathbb{R}$, an ordered knot vector

$$\Xi_k := (0 = \xi_{k,1}, \xi_{k,2}, \ldots, \xi_{k,n_k+p_k+1} = 1), \quad \xi_{k,1} \leq \xi_{k,2} \leq \cdots \leq \xi_{k,n_k+p_k+1},$$

where $\xi_{k,i_k}$ is the $i_k$th knot with $i_k = 1, 2, \ldots, n_k + p_k + 1$ representing the knot index for the coordinate direction $k$. Although not absolutely necessary, assume that $\xi_{k,1} = 0$ and $\xi_{k,n_k+p_k+1} = 1$ for any $k$, so that all parametric intervals are the same as $[0,1]$. The knots may be equally spaced or unequally spaced, resulting in a uniform or non-uniform distribution. More importantly, the knots may be repeated, that is, a knot $\xi_{k,i_k}$ of the knot vector $\Xi_k$ may appear $1 \leq m_{k,i_k} \leq p_k + 1$ times, where $m_{k,i_k}$ is referred to as its multiplicity. The multiplicity has important implications on the regularity properties of B-spline functions, to be discussed later. To monitor knots without repetitions, say, there are $r_k$ distinct knots in $\Xi_k$. Collect them into an auxiliary knot vector $\mathbf{Z}_k := (\zeta_{k,1}, \ldots, \zeta_{k,r_k})$ and define the vector $\mathbf{M}_k := (m_{k,1}, \ldots, m_{k,r_k})$ of their corresponding multiplicities such that

$$\Xi_k = (0 = \underbrace{\zeta_{k,1}, \ldots, \zeta_{k,1}}_{m_{k,1} \text{ times}}, \underbrace{\zeta_{k,2}, \ldots, \zeta_{k,2}}_{m_{k,2} \text{ times}}, \ldots, \underbrace{\zeta_{k,r_k}, \ldots, \zeta_{k,r_k}}_{m_{k,r_k} \text{ times}} = 1), \quad \sum_{i_k=1}^{r_k} m_{k,i_k} = n_k + p_k + 1.$$

A knot vector is called open if its first and last knots appear $p_k + 1$ times. Open knot vectors are standard in CAD [14].

### 2.2. Univariate B-Splines

The B-spline functions for a given degree are defined in a recursive manner using the knot vector. Denote by $N_{i_k,p_k}^k(\xi_k)$ the $i_k$th univariate B-spline function with degree $p_k$ for the coordinate direction $k$. Given the *zero*-degree basis functions,

$$N_{i_k,0}^k(\xi_k) = \begin{cases} 1, & \xi_{k,i_k} \leq \xi_k < \xi_{k,i_k+1}, \\ 0, & \text{otherwise}, \end{cases}$$

for $k = 1, \ldots, d$, all higher-order B-spline functions are efficiently generated by the recursive Cox-de Boor formula [15, 16],

$$N_{i_k,p_k}^k(\xi_k) = \frac{\xi_k - \xi_{k,i_k}}{\xi_{k,i_k+p_k} - \xi_{k,i_k}} N_{i_k,p_k-1}^k(\xi_k) + \frac{\xi_{k,i_k+p_k+1} - \xi_k}{\xi_{k,i_k+p_k+1} - \xi_{k,i_k+1}} N_{i_k+1,p_k-1}^k(\xi_k), \quad (1)$$

where $1 \leq k \leq d$, $1 \leq i_k \leq n_k$, $1 \leq p_k < \infty$, and $0/0$ is considered as *zero*.

The B-spline functions for any $k = 1, \ldots, d$ and $p_k \in \mathbb{N}_0$ satisfy the following desirable properties [13, 14, 15, 16]: (1) they are non-negative, that is, $N_{i_k,p_k}^k(\xi_k) \geq 0$ for all $i_k$ and $\xi_k$; (2) they are locally supported on the interval $[\xi_{k,i_k}, \xi_{k,i_k+p_k+1}]$ for all $i_k$; (3) they are linearly independent, that is, if $\sum_{i_k=1}^{n_k} c^k_{i_k} N^k_{i_k,p_k}(\xi_k) = 0$, then $c^k_{i_k} = 0$ for all $i_k$; (4) they form a partition of unity, that is, $\sum_{i_k=1}^{n_k} N^k_{i_k,p_k}(\xi_k) = 1$, $\xi_k \in [\xi_{k,1}, \xi_{k,n_k+p_k+1}]$; and (5) they are everywhere pointwise $C^\infty$-continuous except at the knots $\xi_{k,i_k}$ of multiplicity $m_{k,i_k}$, where it is $C^{p_k-m_{k,i_k}}$-continuous, provided that $1 \leq m_{k,i_k} < p_k + 1$.

---

[2] The nouns *degree* and *order* associated with IGA are used synonymously in the paper.

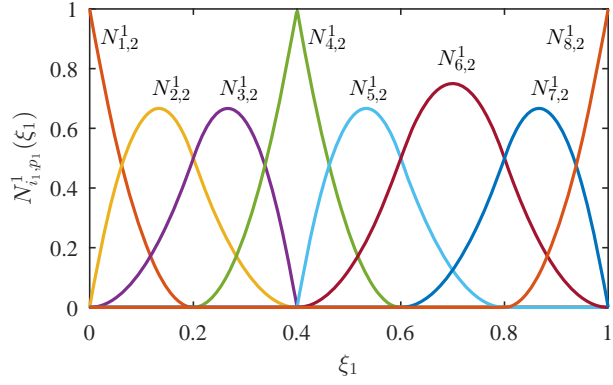

(a)

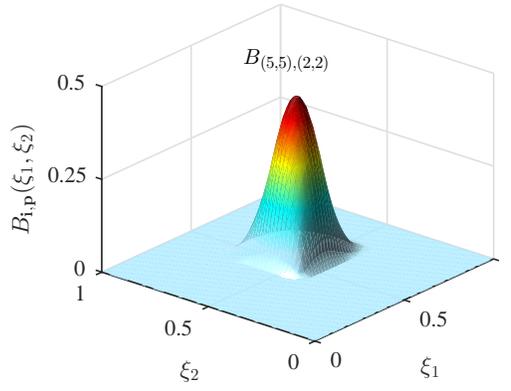

(b)

Figure 1: Quadratic B-splines generated from the knot vectors $\Xi_1 = \Xi_2 = (0, 0, 0, 0.2, 0.4, 0.4, 0.6, 0.8, 1, 1, 1)$ with $n_1 = n_2 = 8$ and $p_1 = p_2 = 2$; (a) eight univariate B-splines for the coordinate direction $\xi_1$; (b) a bivariate B-spline from the tensor product of $N^1_{5,2}(\xi_1)$ and $N^2_{5,2}(\xi_2)$.

Define by
$$\mathcal{N}_k := \mathcal{N}_k(\Xi_k; p_k) := \text{span}\{N^k_{i_k,p_k}(\xi_k)\}_{i_k=1,\ldots,n_k}$$
the space of univariate B-splines with degree $p_k$ for the coordinate direction $k$. Figure 1(a) shows eight univariate quadratic B-spline basis functions $N^1_{i_1,p_1}(\xi_1)$, $i_1 = 1, \ldots, n_1$, when $n_1 = 8$, $p_1 = 2$, and $\Xi_1 = \{0, 0, 0, 0.2, 0.4, 0.4, 0.6, 0.8, 1, 1, 1\}$. The multiplicity of each interior knot is one, except at the fifth knot, where it is two. Therefore, the basis functions are $C^0$-continuous at $\xi_{1,5} = \xi_{1,6} = 0.4$ and $C^1$-continuous at other interior nodes. Clearly, the regularities of B-splines depend on the multiplicities of the knots selected.

*2.3. Multivariate B-Splines*

The multivariate B-splines in $d$ variables $\xi_1, \ldots, \xi_d$ are constructed from the tensor product of the univariate B-splines stemming from the chosen knot vectors $\Xi_1, \ldots, \Xi_d$. The corresponding auxiliary knot vectors and multiplicity vectors are $\mathbf{Z}_1, \ldots, \mathbf{Z}_d$ and $\mathbf{M}_1, \ldots, \mathbf{M}_d$, respectively. A mesh $\mathcal{Q}_h$ in the parametric domain $\hat{\mathcal{D}} = [0,1]^d$ is defined by its partition into $d$-dimensional parametric elements $Q$, that is,

$$\mathcal{Q}_h := \left\{ Q = \otimes_{k=1}^d (\zeta_{k,i_k}, \zeta_{k,i_k+1}) : 1 \leq i_k < r_k - 1 \right\}.$$

Relatedly, if the size of an element $Q \in \mathcal{Q}_h$ is defined as $\hat{h}_Q := \mathrm{diam}(Q)$, then $\hat{h} := \max_{Q \in \mathcal{Q}_h}\{\hat{h}_Q\}$ defines the global mesh size in the parametric domain.

Define two multi-indices $\mathbf{i} := (i_1, \ldots, i_d)$ and $\mathbf{p} := (p_1, \ldots, p_d)$. For the first multi-index, denote by

$$\mathcal{I} := \{ \mathbf{i} = (i_1, \ldots, i_d) : 1 \leq i_k \leq n_k,\ 1 \leq k \leq d \}$$

a multi-index set. Then, for $\mathbf{i} \in \mathcal{I}$ and $\mathbf{p} \in \mathbb{N}_0^d$, the multivariate B-spline function $B_{\mathbf{i},\mathbf{p}} : \hat{\mathcal{D}} \to \mathbb{R}$ is defined as

$$B_{\mathbf{i},\mathbf{p}}(\boldsymbol{\xi}) := \prod_{k=1}^d N_{i_k, p_k}^k(\xi_k) \tag{2}$$

with the corresponding tensor-product B-spline space

$$\mathcal{B}_h := \bigotimes_{k=1}^d \mathcal{N}_k(\boldsymbol{\Xi}_k; p_k) = \bigotimes_{k=1}^d \mathrm{span}\{N_{i_k,p_k}^k(\xi_k)\}_{i_k = 1, \ldots, n_k} = \mathrm{span}\{B_{\mathbf{i},\mathbf{p}}(\boldsymbol{\xi})\}_{\mathbf{i} \in \mathcal{I}}. \tag{3}$$

Note that the functions in $\mathcal{B}_h$ are piecewise polynomials of degree $p_k$ along each coordinate direction $k = 1, \ldots, d$. Figure 1(b) depicts a bivariate quadratic B-spline, which is generated from the knot vectors $\boldsymbol{\Xi}_1 = \boldsymbol{\Xi}_2 = (0, 0, 0, 0.2, 0.4, 0.4, 0.6, 0.8, 1, 1, 1)$ and tensor product of $N_{5,2}^1(\xi_1)$ and $N_{5,2}^2(\xi_2)$.

Due to the tensor-product structure, multivariate B-spline functions inherit most of the aforementioned properties of their univariate counterparts, namely, non-negativity, local support, linear independence, partition of unity, and regularity. The functions are $C^\infty$-continuous in the interior of each element $Q \in \mathcal{Q}_h$, while, across element boundaries, the regularity is decided by the directional regularity in each coordinate.

### 2.4. NURBS

With the multivariate B-spline functions and their space described by (2) and (3), the multivariate NURBS functions and the corresponding space can now be defined using a projective transformation [14, 17]. Associated with each $\mathbf{i} \in \mathcal{I}$, denote by $w_{\mathbf{i}} \in \mathbb{R}^+$ a constant positive weight. As a result, the weight function $w : \hat{\mathcal{D}} \to \mathbb{R}$ can be defined as

$$w(\boldsymbol{\xi}) := \sum_{\mathbf{i} \in \mathcal{I}} w_{\mathbf{i}} B_{\mathbf{i},\mathbf{p}}(\boldsymbol{\xi}).$$

Using the properties of B-splines, it is elementary to show that the weight function is also positive. Given $\mathbf{i} \in \mathcal{I}$ and $\mathbf{p} \in \mathbb{N}_0^d$, the multivariate NURBS function $R_{\mathbf{i},\mathbf{p}} : \hat{\mathcal{D}} \to \mathbb{R}$ is defined as [14, 17]

$$R_{\mathbf{i},\mathbf{p}}(\boldsymbol{\xi}) := \frac{w_{\mathbf{i}} B_{\mathbf{i},\mathbf{p}}(\boldsymbol{\xi})}{w(\boldsymbol{\xi})} = \frac{w_{\mathbf{i}} B_{\mathbf{i},\mathbf{p}}(\boldsymbol{\xi})}{\displaystyle\sum_{\mathbf{i} \in \mathcal{I}} w_{\mathbf{i}} B_{\mathbf{i},\mathbf{p}}(\boldsymbol{\xi})},$$

producing the NURBS function space

$$\mathcal{R}_h := \mathrm{span}\{R_{\mathbf{i},\mathbf{p}}(\boldsymbol{\xi})\}_{\mathbf{i} \in \mathcal{I}} \tag{4}$$

on the parametric domain $\hat{\mathcal{D}}$.

The NURBS functions described in the preceding inherit all of the important properties from their piecewise polynomial counterparts as follows [13]: (1) they constitute a partition of unity; (2) the NURBS and B-splines functions have the same continuity and support; (3) they possess the property of affine transformations; (4) setting all the weights to be equal, a NURBS function reduces to a scaled B-spline function; and (5) the NURBS surfaces and solids are projective transformations of tensor-product, piecewise polynomial entities.

Using multivariate B-splines and NURBS functions, a geometric object in $\mathbb{R}^d$, such as a curve, surface, or solid, can be readily generated. For each $\mathbf{i} \in \mathcal{I}$, let $\mathbf{C_i} \in \mathbb{R}^d$ be a control point. Denote by $n_c := |\mathcal{I}|$ the cardinality of $\mathcal{I}$, representing the number of such control points. Call the collection of such control points $\{\mathbf{C_i}\}_{\mathbf{i} \in \mathcal{I}}$ to be a control mesh. Using NURBS functions, the physical domain $\mathcal{D} \subset \mathbb{R}^d$ is obtained by a geometrical mapping $\mathbf{x} : \hat{\mathcal{D}} \to \mathcal{D} \subset \mathbb{R}^d$, which is described more explicitly by

$$\mathbf{x}(\boldsymbol{\xi}) = \sum_{\mathbf{i} \in \mathcal{I}} R_{\mathbf{i},\mathbf{p}}(\boldsymbol{\xi}) \mathbf{C_i}. \tag{5}$$

A similar mapping can be defined using multivariate B-spline functions. However, not all objects or domains, some of which are commonly used in engineering, can be represented by B-splines. For instance, free-form surfaces and conic sections, such as circles, ellipses, cylinders, spheres, ellipsoids, and tori, cannot be described by piecewise polynomials. In contrast, NURBS functions equipped with judiciously selected weights can represent them exactly [14, 17]. Therefore, the use of NURBS, that is, (5), becomes necessary in the CAD community.

Using the geometrical mapping (5), the physical mesh $\mathcal{K}_h$, say, of the physical domain can now be viewed as the image of the parametric mesh $\mathcal{Q}_h$, that is,

$$\mathcal{K}_h := \{K = \mathbf{x}(Q) : Q \in \mathcal{Q}_h\}, \tag{6}$$

where the element $K$ of the physical mesh is the image of the element $Q$ of the parametric mesh. Following Bazilevs *et el.* [18], define the global physical mesh size as $h := \max_{K \in \mathcal{K}_h} h_K$, where $h_K = \|\boldsymbol{\nabla}\mathbf{x}\|_{L^\infty(Q)} \hat{h}_Q$ is the size of element $K$. Here, $\boldsymbol{\nabla}\mathbf{x}$ is the Jacobian of the mapping $\mathbf{x} : \hat{\mathcal{D}} \to \mathcal{D} \subset \mathbb{R}^d$, comprising a matrix of partial derivatives of the components of $\mathbf{x}$ with respect to those of $\boldsymbol{\xi}$, and $\|\cdot\|_{L^\infty(Q)}$ is the infinity norm of the matrix. Moreover, define the space of NURBS functions in the physical domain $\mathcal{D}$ as the push-forward of the NURBS space $\mathcal{R}_h$ in (4) via

$$\mathcal{V}_h := \operatorname{span}\left\{R_{\mathbf{i},\mathbf{p}} \circ \mathbf{x}^{-1}\right\}_{\mathbf{i} \in \mathcal{I}} = \operatorname{span}\left\{\bar{R}_{\mathbf{i},\mathbf{p}}\right\}_{\mathbf{i} \in \mathcal{I}}, \tag{7}$$

where $\bar{R}_{\mathbf{i},\mathbf{p}} := R_{\mathbf{i},\mathbf{p}} \circ \mathbf{x}^{-1}$ is the NURBS function in the physical domain. It is assumed that the mapping (5) is invertible *almost everywhere* in $\mathcal{D}$ and has smooth inverse on each element $K$ of the physical mesh $\mathcal{K}_h$.

### 2.5. Refinement

The accuracy of IGA depends on the enrichment of the NURBS spaces $\mathcal{R}_h$ and $\mathcal{V}_h$ in (4) and (7) via refinement. There are three principal types of refinement. A simple and straightforward type, namely knot insertion, is equivalent to $h$-refinement commonly used in FEM. For knot insertion, a finer mesh is constructed by adding knots to the existing knot vectors without changing the geometry. As an example, consider inserting a new knot $\xi'_k \in [\xi_{k,l}, \xi_{k,l+1})$, $1 \leq l \leq n_k + p_k$, to the existing knot vector $\boldsymbol{\Xi}_k := (\xi_{k,1}, \xi_{k,2}, \ldots, \xi_{k,n_k+p_k+1})$, which produces $n_k$ B-spline functions. Applying the Cox-de Boor formula (1) to the new knot vector, say,

$$\boldsymbol{\Xi}'_k := (\xi'_{k,1}, \xi'_{k,2}, \ldots, \xi'_{k,n_k+p_k+2}) = (\xi_{k,1}, \xi_{k,2}, \ldots, \xi_{k,l}, \xi'_k, \xi_{k,l+1}, \ldots, \xi_{k,n_k+p_k+1}),$$

a new set of $n_k + 1$ basis functions is created with their span nesting the span of existing basis functions. The process can be repeated for additional knots. Moreover, the $h$-refinement can be performed globally in all $d$ coordinate directions or individually in select coordinate directions.[3] Henceforth, for a NURBS object in $\mathbb{R}^d$, a new set of control points should be defined for the new basis functions to obtain an object that is geometrically and parametrically the same as the original one. In other words, the geometry of the physical domain is preserved. Through knot insertion or $h$-refinement, there are increases in the number of elements as well as in the number of basis functions and, consequently, in the number of control points.

---

[3] The $h$-refinements in all and a single coordinate direction(s) are illustrated in Examples 2 and 3, respectively, of Section 5.

Another prominent type of refinement is called order elevation, which is reminiscent of $p$-refinement in FEM. For order elevation, higher-order polynomials are used on the same mesh, again, without changing the geometry. In other words, the geometry of the physical domain is also upheld in $p$-refinement. A third type of refinements, called $k$ refinement, has no analog in FEM, but it encompasses both knot insertion and order elevation. For a detailed description of all three types of refinement, including examples, read the book by Cottrell *et al.* [12].

2.6. Error Estimate of NURBS

Let $L^2(\mathcal{D})$ be a Hilbert space of real-valued square-integrable functions on $\mathcal{D} \subset \mathbb{R}^d$ equipped with the usual inner product and induced norm

$$(u,v)_{L^2(\mathcal{D})} := \int_{\mathcal{D}} u(\mathbf{x})v(\mathbf{x})d\mathbf{x} \quad \text{and} \quad \|u\|_{L^2(\mathcal{D})} := (u,u)_{L^2(\mathcal{D})}^{1/2},$$

respectively, for any $u, v \in L^2(\mathcal{D})$. Define a Sobolev space $H^l(\mathcal{D})$ of order $l \in \mathbb{N}_0$, where $H^0(\mathcal{D}) = L^2(\mathcal{D})$, with the standard norm and seminorm

$$\|u\|_{H^l(\mathcal{D})} := \left(\sum_{k=0}^{l} |u|_{H^k(\mathcal{D})}^2\right)^{1/2} \quad \text{and} \quad |u|_{H^k(\mathcal{D})} := \left(\sum_{|\boldsymbol{i}|=k} \|\partial^{\boldsymbol{i}} u\|_{L^2(\mathcal{D})}^2\right)^{1/2},$$

respectively, for any $u \in H^l(\mathcal{D})$, where the multi-index $\boldsymbol{i} = (i_1, \ldots, i_d) \in \mathbb{N}^d$, $|\boldsymbol{i}| = i_1 + \cdots + i_d$, and $\partial^{\boldsymbol{i}} u = \partial^{i_1 + \cdots + i_d} / \partial x_1^{i_1} \cdots x_d^{i_d}$.

For error analysis of NURBS-based IGA, consider interpolating a function $u \in L^2(\mathcal{D})$ defined on $\mathcal{D} \in \mathbb{R}^d$. However, the decay of interpolation error depends on how the NURBS space is refined. Here, an important result pertaining to the error analysis of NURBS, derived by Bazilevs *et al.* [18], is briefly summarized. In doing so, consider a family of meshes $\{\mathcal{Q}_h\}_{h>0}$ over the parametric domain $\tilde{\mathcal{D}}$, which is subjected to $h$-refinement, while keeping the polynomial degree fixed. The error estimate is based on introducing a support extension $\bar{Q}$ of an element $Q$ of the parametric mesh $\mathcal{Q}_h$ defined as the union of the supports of basis functions whose supports intersect the element $Q$. Similarly, the physical support extension $\bar{K}$, say, of an element $K = \mathbf{x}(Q)$ of the physical mesh $\mathcal{K}_h$ in (6) is obtained as the image of $Q$ via the geometrical mapping in (5). Given a function $u \in L^2(\mathcal{D})$, denote by $\Pi_{\mathcal{V}_h} : L^2(\mathcal{D}) \to \mathcal{V}_h$ a projective operator, where the NURBS space $\mathcal{V}_h$ is defined in (7). The local error, described in Theorem 1, proves convergence of NURBS-based function interpolations.

**Theorem 1** (Local Error Estimate [18]). *Let $\mathcal{V}_h$ be the NURBS space in the physical domain $\mathcal{D}$ defined by the NURBS functions endowed with degrees $p_1 = \cdots = p_d = p$, where $p \in \mathbb{N}_0$. Given the integers $l$ and $r$ such that $0 \leq l \leq r \leq p+1$, the local estimate of the interpolation error for a function $u \in L^2(\mathcal{D}) \cap H^r(\bar{K})$, measured in terms of the $l$th-order seminorm $|\cdot|_{H^l(K)}$ of the Sobolev space $H^l(K)$, is*

$$|u - \Pi_{\mathcal{V}_h} u|_{H^l(K)} \leq C \sum_{K \in \mathcal{K}_h} h_K^{r-l} \sum_{i=0}^{r} \|\boldsymbol{\nabla} \mathbf{x}\|_{L^\infty(\bar{Q})}^{i-r} |u|_{H^i(\bar{K})}^2,$$

*where $h_K$ is the size of element $K$ of $\mathcal{K}_h$ and $C$ is a constant that depends on $p$ and the shape of the domain $\mathcal{D}$.*

In addition, a follow-up global error estimate, described by Theorem 2 of the paper by Bazilevs *et al.* [18], is given. These error estimates demonstrate that the NURBS space delivers the optimal rate of convergence, as for the finite-element spaces of the same degree.

# 3. Karhunen-Loève Representation

Let $(\Omega, \mathcal{F}, \mathbb{P})$ be a complete probability space, where $\Omega$ is a sample space, $\mathcal{F}$ is a $\sigma$-field on $\Omega$, and $P : \mathcal{F} \to [0,1]$ is a probability measure. Denote by $L^2(\Omega, \mathcal{F}, \mathbb{P})$ a Hilbert space of random variables defined on $(\Omega, \mathcal{F}, \mathbb{P})$ and by $L^2(\mathcal{D} \times \Omega)$ a Hilbert space of random fields defined on $\mathcal{D}$. A random variable or random field, if it is a member of the associated Hilbert space, has finite second-moment properties.

## 3.1. Random Field

A real-valued random field $\alpha$ defined on a bounded domain $\mathcal{D} \subset \mathbb{R}^d$, where $d = 1$, 2, or 3, is a mapping $\alpha : \mathcal{D} \times \Omega \to \mathbb{R}$ such that for each $\mathbf{x} \in \mathcal{D}$, $\alpha(\mathbf{x}, \cdot)$ is a random variable with respect to $(\Omega, \mathcal{F}, \mathbb{P})$. Given the expectation operator $\mathbb{E}$ with respect to the probability measure $\mathbb{P}$, denote by $\mu(\mathbf{x}) := \mathbb{E}[\alpha(\mathbf{x}, \cdot)]$ the mean function and by

$$\Gamma(\mathbf{x}, \mathbf{x}') := \mathbb{E}[(\alpha(\mathbf{x}, \cdot) - \mu(\mathbf{x}))(\alpha(\mathbf{x}', \cdot) - \mu(\mathbf{x}'))], \quad \mathbf{x}, \mathbf{x}' \in \mathcal{D},$$

the covariance function of $\alpha(\mathbf{x}, \cdot)$. Without loss of generality, assume that $\mu(\mathbf{x}) = 0$. A non-*zero*-mean random field can be obtained by just adding the mean function to a *zero*-mean random field.

The random fields are often assumed to be homogeneous or stationary, meaning that their finite-dimensional probability distributions are invariant under arbitrary translations. This implies that the covariance function is a function of the argument difference $\mathbf{x} - \mathbf{x}'$. Moreover, random fields are sometimes assumed to be isotropic, that is, invariant under orthogonal transformations. In which case, the covariance function becomes a function of the distance $\|\mathbf{x} - \mathbf{x}'\|$. Finally, additional assumptions are needed to ensure that the samples of random fields are continuous and differentiable in a mean-square or almost sure sense. Figure 2 shows four commonly used covariance functions of a homogeneous random field defined on $[0, 1]$, comprising exponential, Gaussian, sinusoidal, and Bessel functions. Each covariance function contains parameters describing the variance $\sigma^2$ and correlation length parameter $b$ of the random field.

## 3.2. Karhunen-Loève Expansion

Let $\alpha(\mathbf{x}, \cdot) \in L^2(\mathcal{D} \times \Omega)$ be a random field such that for each $\mathbf{x} \in \mathcal{D}$, $\alpha(\mathbf{x}, \cdot)$ is a random variable in $L^2(\Omega, \mathcal{F}, \mathbb{P})$ and, given a realization $\omega \in \Omega$, $\alpha(\mathbf{x}, \omega) \in L^2(\mathcal{D})$. A definition of the Hilbert-Schmidt integral operator [19], followed by a few relevant propositions and a theorem, leads to a formal description of the KL expansion.

**Definition 2.** *For a bounded domain $\mathcal{D} \in \mathbb{R}^d$, a function $\kappa : \mathcal{D} \times \mathcal{D} \to \mathbb{R}$ is called a Hilbert-Schmidt kernel if $\kappa \in L^2(\mathcal{D} \times \mathcal{D})$, that is, if*

$$\int_{\mathcal{D}} \int_{\mathcal{D}} |\kappa(\mathbf{x}, \mathbf{x}')|^2 d\mathbf{x} d\mathbf{x}' < \infty.$$

*Correspondingly, the Hilbert-Schmidt operator $\mathcal{G}_\kappa : L^2(\mathcal{D}) \to L^2(\mathcal{D})$ is defined as*

$$(\mathcal{G}_\kappa \phi)(\mathbf{x}) := \int_{\mathcal{D}} \kappa(\mathbf{x}, \mathbf{x}') \phi(\mathbf{x}') d\mathbf{x}' \quad \forall \ \phi(\mathbf{x}) \in L^2(\mathcal{D}).$$

**Proposition 3.** *The covariance function $\Gamma : \mathcal{D} \times \mathcal{D} \to \mathbb{R}$ of a random field $\alpha(\mathbf{x}, \cdot) \in L^2(\mathcal{D} \times \Omega)$ is a Hilbert-Schmidt kernel, and, therefore, $\mathcal{G}_\Gamma : L^2(\mathcal{D}) \to L^2(\mathcal{D})$, defined as*

$$(\mathcal{G}_\Gamma \phi)(\mathbf{x}) := \int_{\mathcal{D}} \Gamma(\mathbf{x}, \mathbf{x}') \phi(\mathbf{x}') d\mathbf{x}' \quad \forall \ \phi(\mathbf{x}) \in L^2(\mathcal{D}), \tag{8}$$

*is the Hilbert-Schmidt operator associated with the covariance function.*

*Proof.* Since $\alpha(\mathbf{x}, \cdot) \in L^2(\mathcal{D} \times \Omega)$, the covariance function is square-integrable and hence a bounded function. Therefore, the results of the proposition follow readily. □

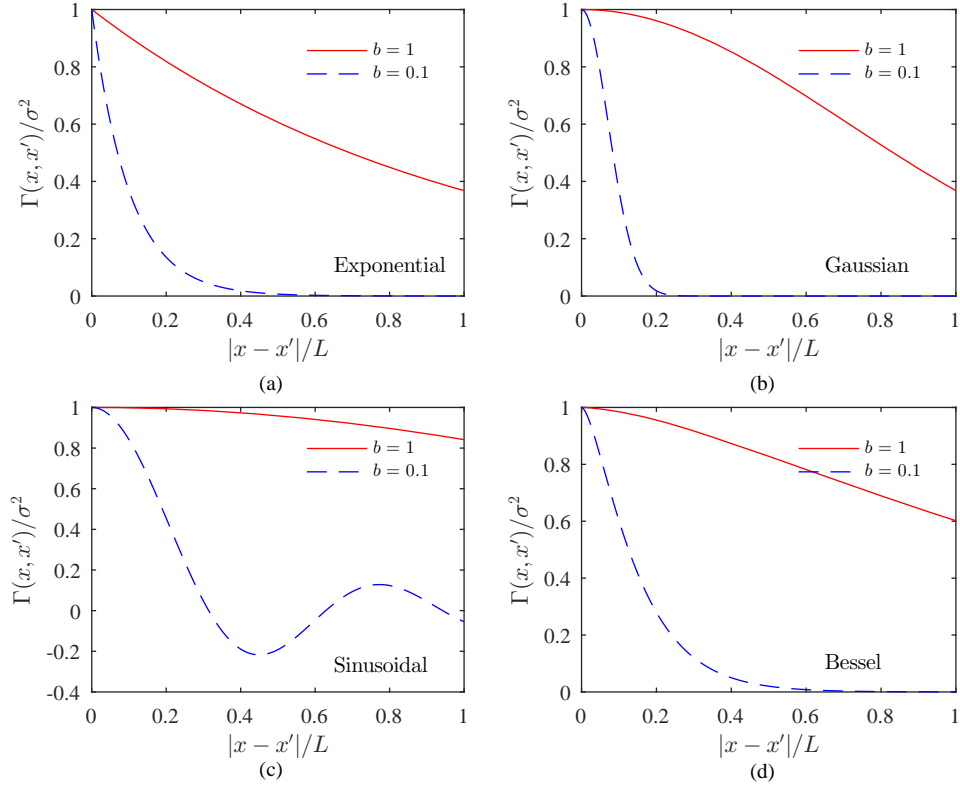

Figure 2: Four commonly used covariance functions of a homogeneous random field defined on $[0, 1]$; (a) exponential function: $\Gamma(x, x') = \sigma^2 \exp[-|x - x'|/(bL)]$; (b) Gaussian function: $\Gamma(x, x') = \sigma^2 \exp[-|x - x'|^2/(bL)^2]$; and (c) sinusoidal function: $\Gamma(x, x') = \sigma^2/[|x - x'|/(bL)] \sin(|x - x'|/(bL))$; and (d) Bessel function: $\Gamma(x, x') = \sigma^2 |x - x'|/(bL) K_1(|x - x'|/(bL))$. Here, $L = 1$; $K_1$ is the modified second-kind Bessel function of order one.

**Proposition 4.** *Let $\mathcal{G}_\Gamma : L^2(\mathcal{D}) \to L^2(\mathcal{D})$ defined in (8) be the Hilbert-Schmidt operator associated with the covariance function $\Gamma : \mathcal{D} \times \mathcal{D} \to \mathbb{R}$ of a zero-mean random field $\alpha(\mathbf{x}, \cdot) \in L^2(\mathcal{D} \times \Omega)$. Then $\mathcal{G}_\Gamma$ is a linear, compact, positive-semidefinite, and self-adjoint operator.*

*Proof.* From (8), $\mathcal{G}_\Gamma$ is obviously linear. Any Hilbert-Schmidt operator from $L^2(\mathcal{D})$ to $L^2(\mathcal{D})$ can be expressed as the limit of a sequence of bounded finite-rank operators. Then, from Lemma 1.2.3 of Atkinson [20], $\mathcal{G}_\Gamma$ is compact. To prove positive-semidefiniteness, invoke the usual inner product $(\cdot, \cdot)_{L^2(\mathcal{D})}$ of $L^2(\mathcal{D})$ to demonstrate, for any $0 \neq \phi(\mathbf{x}) \in L^2(\mathcal{D})$, that

$$\begin{aligned}
((\mathcal{G}_\Gamma \phi)(\mathbf{x}), \phi(\mathbf{x}))_{L^2(\mathcal{D})} &= \int_\mathcal{D} \int_\mathcal{D} \Gamma(\mathbf{x}, \mathbf{x}') \phi(\mathbf{x}) \phi(\mathbf{x}') d\mathbf{x} d\mathbf{x}' \\
&= \int_\mathcal{D} \int_\mathcal{D} \mathbb{E}\left[\alpha(\mathbf{x}, \cdot) \alpha(\mathbf{x}', \cdot)\right] \phi(\mathbf{x}) \phi(\mathbf{x}') d\mathbf{x} d\mathbf{x}' \\
&= \mathbb{E}\left[\left(\int_\mathcal{D} \alpha(\mathbf{x}, \cdot) \phi(\mathbf{x}) d\mathbf{x}\right) \left(\int_\mathcal{D} \alpha(\mathbf{x}', \cdot) \phi(\mathbf{x}') d\mathbf{x}'\right)\right] \\
&= \mathbb{E}\left[\left(\int_\mathcal{D} \alpha(\mathbf{x}, \cdot) \phi(\mathbf{x}) d\mathbf{x}\right)^2\right] \\
&\geq 0,
\end{aligned}$$

where Fubini's theorem is employed to interchange the integrals.

Finally, as the covariance function $\Gamma(\mathbf{x}, \mathbf{x}')$ is symmetric with respect to arguments $\mathbf{x}$ and $\mathbf{x}'$,

$$((\mathcal{G}_\Gamma \phi)(\mathbf{x}), \psi(\mathbf{x}))_{L^2(\mathcal{D})} = (\phi(\mathbf{x}), (\mathcal{G}_\Gamma \psi)(\mathbf{x}))_{L^2(\mathcal{D})}$$

for any $\phi(\mathbf{x}), \psi(\mathbf{x}) \in L^2(\mathcal{D})$, proving that $\mathcal{G}_\Gamma$ is self-adjoint. $\square$

A linear, compact, positive-semidefinite, and self-adjoint operator, such as $\mathcal{G}_\Gamma : L^2(\mathcal{D}) \to L^2(\mathcal{D})$, on an infinite dimensional Hilbert space has spectral properties resembling those of a positive-semidefinite, symmetric matrix. Proposition 5 describes the spectral properties of $\mathcal{G}_\Gamma$.

**Proposition 5.** *Let $\mathcal{G}_\Gamma : L^2(\mathcal{D}) \to L^2(\mathcal{D})$ defined in (8) be the Hilbert-Schmidt operator associated with the covariance function $\Gamma : \mathcal{D} \times \mathcal{D} \to \mathbb{R}$ of a random field $\alpha(\mathbf{x}, \cdot) \in L^2(\mathcal{D} \times \Omega)$. There exists an infinite sequence of eigenpairs $\{\lambda_i, \phi_i(\mathbf{x})\}_{i \in \mathbb{N}}$ of $\mathcal{G}_\Gamma$, which is the solution of*

$$(\mathcal{G}_\Gamma \phi)(\mathbf{x}) = \lambda \phi(\mathbf{x}) \quad or \quad \int_\mathcal{D} \Gamma(\mathbf{x}, \mathbf{x}') \phi(\mathbf{x}') d\mathbf{x}' = \lambda \phi(\mathbf{x}), \tag{9}$$

*known as the Fredholm integral equation of the second kind. Moreover, the eigensolutions, where the eigenfunctions have been normalized such that $\|\phi_i(\mathbf{x})\|^2_{L^2(\mathcal{D})} := \int_\mathcal{D} \phi_i^2(\mathbf{x}) d\mathbf{x} = 1$, satisfy the following properties.*

1. *The eigenvalues $\lambda_i \in \mathbb{R}_0^+$, $i \in \mathbb{N}$, are real and non-negative having zero as the only point of accumulation. Moreover, the eigenvalues can be arranged in a descending order as follows: $\lambda_1 \geq \lambda_2 \geq \cdots \geq 0$.*

2. *The eigenfunctions $\phi_i(\mathbf{x}) \in L^2(\mathcal{D})$, $i \in \mathbb{N}$, corresponding to distinct eigenvalues, are mutually orthonormal, that is,*

$$(\phi_i(\mathbf{x}), \phi_j(\mathbf{x}))_{L^2(\mathcal{D})} := \int_\mathcal{D} \phi_i(\mathbf{x}) \phi_j(\mathbf{x}) d\mathbf{x} = \delta_{ij}, \tag{10}$$

   *where $\delta_{ij}$ is the Kronecker delta, that is, $\delta_{ij} = 1$ when $i = j$ and $\delta_{ij} = 0$ when $i \neq j$.*

3. *The number of eigenfunctions corresponding to non-zero eigenvalues is finite.*

4. *The sequence of eigenfunctions $\{\phi_i(\mathbf{x})\}_{i \in \mathbb{N}}$ forms an orthonormal basis of $L^2(\mathcal{D})$, that is, $L^2(\mathcal{D}) = \mathrm{span}\{\phi_i(\mathbf{x})\}_{i \in \mathbb{N}}$.*

**Theorem 6** (Mercer's theorem [21]). *Let $\Gamma : \mathcal{D} \times \mathcal{D} \to \mathbb{R}$ be a continuous covariance function of a random field $\alpha(\mathbf{x}, \cdot) \in L^2(\mathcal{D} \times \Omega)$ and $\mathcal{G}_\Gamma : L^2(\mathcal{D}) \to L^2(\mathcal{D})$ defined in (8) be the associated Hilbert-Schmidt operator. If $\{\lambda_i, \phi_i(\mathbf{x})\}_{i \in \mathbb{N}}$ is an infinite sequence of eigenpairs of $\mathcal{G}_\Gamma$, then*

$$\Gamma(\mathbf{x}, \mathbf{x}') = \sum_{i=1}^\infty \lambda_i \phi_i(\mathbf{x}) \phi_i(\mathbf{x}') \tag{11}$$

*for all $\mathbf{x}, \mathbf{x}' \in \mathcal{D}$, where the infinite series on the right converges pointwise and uniformly on $\mathcal{D} \times \mathcal{D}$.*

Given the mathematical results of Propositions 3, 4, and 5 and Theorem 6, the KL expansion is presented as follows.

**Theorem 7** (Karhunen-Loève). *Let $\alpha(\mathbf{x}, \cdot) \in L^2(\mathcal{D} \times \Omega)$ be a real-valued random field with zero mean, continuous covariance function $\Gamma : \mathcal{D} \times \mathcal{D} \to \mathbb{R}$, and associated Hilbert-Schmidt operator $\mathcal{G}_\Gamma : L^2(\mathcal{D}) \to L^2(\mathcal{D})$ defined in (8). Given an infinite sequence of eigenpairs $\{\lambda_i, \phi_i(\mathbf{x})\}_{i \in \mathbb{N}}$ of $\mathcal{G}_\Gamma$, the random field admits a convergent infinite series expansion*

$$\alpha(\mathbf{x}, \cdot) \sim \sum_{i=1}^\infty \sqrt{\lambda_i} \phi_i(\mathbf{x}) X_i, \tag{12}$$

*where $\{X_i\}_{i \in \mathbb{N}}$ is an infinite sequence of zero-mean, standardized, uncorrelated random variables, that is,*

$$\mathbb{E}[X_i] = \int_\Omega X_i(\omega) d\mathbb{P}(\omega) = 0,$$

$$\mathbb{E}[X_i X_j] = \int_\Omega X_i(\omega) X_j(\omega) d\mathbb{P}(\omega) = \delta_{ij} \quad i, j \in \mathbb{N},$$

with each random variable $X_i$ defined, for $\lambda_i \neq 0$, as

$$X_i := \frac{1}{\sqrt{\lambda_i}} \int_\mathcal{D} \alpha(\mathbf{x}, \cdot) \phi_i(\mathbf{x}) d\mathbf{x}, \quad \mathbf{x} \in \mathcal{D}. \tag{13}$$

*Proof.* With the recognition that $\alpha(\mathbf{x}, \cdot)$ has *zero* mean, apply the expectation operator on (13) to obtain $\mathbb{E}[X_i] = 0$ for all $i$. The orthogonality of the random variables $\{X_i\}_{i \in \mathbb{N}}$ follows from the orthogonality of the eigenfunctions of $\mathcal{G}_\Gamma$. Indeed, using (13) again, with (9) and (10) in mind,

$$\begin{aligned}
\mathbb{E}[X_i X_j] &= \frac{1}{\sqrt{\lambda_i}} \frac{1}{\sqrt{\lambda_j}} \mathbb{E}\left[\int_\mathcal{D} \int_\mathcal{D} \alpha(\mathbf{x}, \cdot) \alpha(\mathbf{x}', \cdot) \phi_i(\mathbf{x}) \phi_j(\mathbf{x}') d\mathbf{x} d\mathbf{x}'\right] \\
&= \frac{1}{\sqrt{\lambda_i}} \frac{1}{\sqrt{\lambda_j}} \int_\mathcal{D} \left(\int_\mathcal{D} \Gamma(\mathbf{x}, \mathbf{x}') \phi_j(\mathbf{x}') d\mathbf{x}'\right) \phi_i(\mathbf{x}) d\mathbf{x} \\
&= \frac{1}{\sqrt{\lambda_i}} \frac{1}{\sqrt{\lambda_j}} \int_\mathcal{D} \lambda_j \phi_j(\mathbf{x}) \phi_i(\mathbf{x}) d\mathbf{x} \\
&= \delta_{ij}.
\end{aligned}$$

For convergence analysis, given an integer $N \in \mathbb{N}$, define a second-moment error

$$e_N(\mathbf{x}) := \mathbb{E}\left[\left\{\alpha(\mathbf{x}, \cdot) - \sum_{i=1}^N \sqrt{\lambda_i} \phi_i(\mathbf{x}) X_i\right\}^2\right]$$

committed by the $N$-term truncation of the infinite series on the right side of (12). On expansion,

$$\begin{aligned}
e_N(\mathbf{x}) &= \mathbb{E}\left[\alpha^2(\mathbf{x}, \cdot)\right] + \mathbb{E}\left[\left\{\sum_{i=1}^N \sqrt{\lambda_i} \phi_i(\mathbf{x}) X_i\right\}^2\right] - 2\mathbb{E}\left[\alpha(\mathbf{x}, \cdot) \sum_{i=1}^N \sqrt{\lambda_i} \phi_i(\mathbf{x}) X_i\right] \\
&= \Gamma(\mathbf{x}, \mathbf{x}) + \mathbb{E}\left[\sum_{i=1}^N \sum_{j=1}^N \sqrt{\lambda_i} \sqrt{\lambda_j} \phi_i(\mathbf{x}) \phi_j(\mathbf{x}) X_i X_j\right] - 2\mathbb{E}\left[\sum_{i=1}^N \int_\mathcal{D} \alpha(\mathbf{x}, \cdot) \alpha(\mathbf{x}', \cdot) \phi_i(\mathbf{x}') \phi_i(\mathbf{x}) d\mathbf{x}'\right] \\
&= \Gamma(\mathbf{x}, \mathbf{x}) + \sum_{i=1}^N \lambda_i \phi_i^2(\mathbf{x}) - 2\left[\sum_{i=1}^N \int_\mathcal{D} \mathbb{E}[\alpha(\mathbf{x}, \cdot) \alpha(\mathbf{x}', \cdot)] \phi_i(\mathbf{x}') \phi_i(\mathbf{x}) d\mathbf{x}'\right] \\
&= \Gamma(\mathbf{x}, \mathbf{x}) + \sum_{i=1}^N \lambda_i \phi_i^2(\mathbf{x}) - 2\left[\sum_{i=1}^N \left(\int_\mathcal{D} \Gamma(\mathbf{x}, \mathbf{x}') \phi_i(\mathbf{x}') d\mathbf{x}'\right) \phi_i(\mathbf{x})\right] \\
&= \Gamma(\mathbf{x}, \mathbf{x}) + \sum_{i=1}^N \lambda_i \phi_i^2(\mathbf{x}) - 2 \sum_{i=1}^N \lambda_i \phi_i^2(\mathbf{x}) \\
&= \Gamma(\mathbf{x}, \mathbf{x}) - \sum_{i=1}^N \lambda_i \phi_i^2(\mathbf{x}).
\end{aligned}$$

Here, the second line is obtained by recognizing that the variance of $\alpha(\mathbf{x}, \cdot)$ is $\Gamma(\mathbf{x}, \mathbf{x})$ and applying (13); the third line is attained by using orthogonality of random variables and bringing the expectation operator inside the integral; the fourth and fifth lines are acquired by definition of covariance function and applying (9) and (10); and finally, the last line is the result of reduction. Taking the limit $N \to \infty$ and invoking Mercer's theorem yields

$$\lim_{N \to \infty} e_N(\mathbf{x}) = \Gamma(\mathbf{x}, \mathbf{x}) - \sum_{i=1}^\infty \lambda_i \phi_i^2(\mathbf{x}) = 0,$$

proving that the right side of (12) converges to $\alpha(\mathbf{x},\cdot)$ uniformly on $\mathcal{D}$ and in $L^2(\Omega,\mathcal{F},\mathbb{P})$. □

*3.3. Karhunen-Loève Approximation*

The KL expansion contains an infinite number of eigenpairs and random variables. In practice, the number must be finite, meaning that the expansion must be truncated. A straightforward approach, assuming that the eigenvalues have been arranged in a descending sequence, entails retaining the first $N \in \mathbb{N}$ terms of the expansion. The result is an $N$-term truncation or KL approximation

$$\alpha_N(\mathbf{x},\cdot) = \sum_{i=1}^{N} \sqrt{\lambda_i}\phi_i(\mathbf{x})X_i \qquad (14)$$

of $\alpha(\mathbf{x},\cdot)$, comprising eigenpairs $\{\lambda_i, \phi_i(\mathbf{x})\}_{1 \leq i \leq N}$ and random variables $X_i$, $i = 1, \ldots, N$. This is commonly referred to as random field discretization. In consequence, the statistical variation of random field $\alpha(\mathbf{x},\cdot)$ is being swapped with those possessed by $N$ uncorrelated random variables $X_1, \ldots, X_N$. Therefore, the value of $N$ should be selected judiciously not only for maintaining desired accuracy in the discretization, but also for computational expediency.

To determine the quality of a KL approximation, a simple and efficient approach entails second-moment analysis. As already alluded to in the proof of Theorem 7, the second-moment error of $\alpha_N(\mathbf{x},\cdot)$ can be expressed by

$$\bar{e}_N(\mathbf{x}) := \frac{e_N(\mathbf{x})}{\sigma^2(\mathbf{x})} = 1 - \frac{1}{\sigma^2(\mathbf{x})} \sum_{i=1}^{N} \lambda_i \phi_i^2(\mathbf{x}) \qquad (15)$$

when normalized with respect to the variance $\sigma^2(\mathbf{x}) = \Gamma(\mathbf{x},\mathbf{x})$ of $\alpha(\mathbf{x},\cdot)$. However, the error measure in (15) is given at a specific point $\mathbf{x} \in \mathcal{D}$ and is, therefore, local. Henceforth, a global error measure, obtained by averaging over all points of $\mathcal{D}$, can also be obtained as

$$\tilde{e}_N := \frac{1}{|\mathcal{D}|} \int_{\mathcal{D}} \bar{e}_N(\mathbf{x}) d\mathbf{x} = 1 - \frac{1}{|\mathcal{D}|} \sum_{i=1}^{N} \lambda_i \int_{\mathcal{D}} \frac{\phi_i^2(\mathbf{x})}{\sigma^2(\mathbf{x})} d\mathbf{x}, \qquad (16)$$

where $|\mathcal{D}| := \int_{\mathcal{D}} d\mathbf{x}$ is the Lebesgue measure, such as length for $d = 1$, area for $d = 2$, or volume for $d = 3$, of $\mathcal{D}$. If, in addition, the variance is constant, say, $\sigma^2$, as is the case for homogeneous random fields, then (16) reduces to

$$\tilde{e}_N := 1 - \frac{1}{|\mathcal{D}|}\frac{1}{\sigma^2} \sum_{i=1}^{N} \lambda_i. \qquad (17)$$

Both the local and global error measures in (15)-(17) indicate that the truncated KL expansion always underestimates the variance of the original random field. Clearly, the variance of the KL approximation $\alpha_N(\mathbf{x},\cdot)$ is obtained by adding individual contributions from all $N$ eigenmodes. Therefore, the larger the value of $N$, the smaller the respective error. More importantly, given a value of $N$, the effectiveness of an $N$-term KL approximation is predicated on how fast the eigenvalues decay with respect to the eigenmodes. The rate of decay depends strongly on the properties of the covariance function, especially the correlation length parameter of the covariance function. The smoothness of the covariance function also determines the rate of eigenvalue decay and regularity of eigenfunctions [22]. Nonetheless, a remarkable property of the KL expansion is its error-minimizing property; that is, given a fixed $N$, the KL approximation in (14) has been proven to be optimal among all series expansion methods with respect to a global mean-square error [7].

*3.4. Remarks*

The KL expansion is useful for a number of reasons. First, the approximation holds for both homogeneous and inhomogeneous fields. Second, it is optimal in the sense that the mean-square error of the approximation is minimized. Third, the sequence of KL approximations $\{\alpha_N(\mathbf{x},\cdot)\}_{N \in \mathbb{N}}$ of $\alpha(\mathbf{x},\cdot)$ converges in mean-square to the correct limit as $N \to \infty$ at each $\mathbf{x} \in \mathcal{D}$. Finally, as $N \to \infty$, the covariance function of $\alpha_N$ approaches the covariance function of $\alpha$, rendering $\alpha_N$ equal to $\alpha$ in the second-moment sense. If $\alpha$ is a Gaussian random field, then $X_i$, $i = 1, \ldots, N$, in (14) are independent standard Gaussian random variables. In consequence, each element of the sequence $\{\alpha_N\}_{N \in \mathbb{N}}$ is a Gaussian random field.

However, if $\alpha$ is a non-Gaussian random field, then $X_i$, $i = 1, \ldots, N$, are uncorrelated yet dependent non-Gaussian random variables. In which case, finding their probability distribution is not trivial. One class of non-Gaussian random fields for which the use of KL expansion can possibly be exploited is the class of translation random fields, where a non-Gaussian random field is defined as a nonlinear, memoryless transformation of a Gaussian random field [23]. But, even then, there are conditions on the covariance properties that must be fulfilled before proceeding with the transformation [11, 23]. Finally, it should be noted that the KL expansion uses only the second-moment properties, such as the covariance function, of a random field. If two non-Gaussian random fields have identical covariance functions, but significantly different higher-order moments or higher-order finite-dimensional distributions, then they will have the same KL representation, but their sample properties may vary markedly. This has important ramifications in reliability analysis.

## 4. Galerkin Isogeometric Method

The implementation of a KL approximation is predicated on the knowledge of eigensolutions of the integral eigenvalue problem defined by (9). However, analytical or exact solutions of the eigenvalue problem exist when the covariance function $\Gamma : \mathcal{D} \times \mathcal{D} \to \mathbb{R}$ is separable and has simpler functional forms, such as exponential functions, or the domain $\mathcal{D}$ is rectangular. For arbitrary covariance functions or arbitrary domains, numerical methods are often needed to solve the eigenvalue problem. In this section, the basis functions from isogeometric analysis, such as B-splines and NURBS, in conjunction with Galerkin's finite-dimensional approximation, are exploited to solve the eigenvalue problem.

### 4.1. Galerkin Discretization

The variational or weak formulation for solving the Fredholm integral equation (9) entails finding an eigenpair $\{\lambda, \phi(\mathbf{x})\} \subset \mathbb{R}_0^+ \times L^2(\mathcal{D})$ such that

$$((\mathcal{G}_\Gamma \phi)(\mathbf{x}), \psi(\mathbf{x}))_{L^2(\mathcal{D})} = \lambda (\phi(\mathbf{x}), \psi(\mathbf{x}))_{L^2(\mathcal{D})} \tag{18}$$

for any $\psi(\mathbf{x}) \in L^2(\mathcal{D})$. However, since (9) is, in general, not exactly solvable, neither is (18). Therefore, (18) must be solved approximately, say, by Galerkin discretization.

Let $\{\mathcal{S}_h\}_{h>0}$ be a sequence of finite-dimensional approximating subspaces of $L^2(\mathcal{D})$. Denote by $P_h$ a projection of $L^2(\mathcal{D})$ onto $\mathcal{S}_h$. Then, a Galerkin solution calls for finding an eigenpair $\{\lambda_h, \phi_h(\mathbf{x})\} \subset \mathbb{R}_0^+ \times \mathcal{S}_h$ such that

$$((\mathcal{G}_\Gamma \phi_h)(\mathbf{x}), \psi(\mathbf{x}))_{L^2(\mathcal{D})} = (\lambda_h \phi_h(\mathbf{x}), \psi(\mathbf{x}))_{L^2(\mathcal{D})} \tag{19}$$

for any $\psi(\mathbf{x}) \in \mathcal{S}_h \subset L^2(\mathcal{D})$ or, equivalently, solving

$$\int_\mathcal{D} \left( \int_\mathcal{D} \Gamma(\mathbf{x}, \mathbf{x}') \phi_h(\mathbf{x}') d\mathbf{x}' \right) \psi(\mathbf{x}) d\mathbf{x} = \lambda_h \int_\mathcal{D} \phi_h(\mathbf{x}) \psi(\mathbf{x}) d\mathbf{x}, \ \forall \ \psi(\mathbf{x}) \in \mathcal{S}_h \subset L^2(\mathcal{D}). \tag{20}$$

This is known as the Galerkin variational or weak form of (9) relative to the subspace $\mathcal{S}_h$. In general, the solution of (19) or (20) is an approximate solution of (9). The existence and uniqueness of the Galerkin solution of the integral equations of the second kind were discussed by Atkinson [20].

A standard error analysis of $\lambda_h$ and $\phi_h(\mathbf{x})$ entails demonstrating that [20]

$$\|\mathcal{G}_\Gamma - P_h \mathcal{G}_\Gamma\|_{L^2(\mathcal{D})} \to 0 \text{ as } h \to 0 \tag{21}$$

with respect to the norm of $L^2(\mathcal{D})$. For a general $\phi(\mathbf{x}) \in L^2(\mathcal{D})$, define $P_h \phi(\mathbf{x})$ to be the solution of the following minimization problem:

$$\|\phi(\mathbf{x}) - P_h\phi(\mathbf{x})\|_{L^2(\mathcal{D})} := \min_{\psi(\mathbf{x}) \in \mathcal{S}_h} \|\phi(\mathbf{x}) - \psi(\mathbf{x})\|_{L^2(\mathcal{D})}. \qquad (22)$$

Since $\mathcal{S}_h$ is finite-dimensional, (22) has a solution. Moreover, by $\mathcal{S}_h$ being an inner product space, the solution is unique. If, indeed,

$$P_h\phi(\mathbf{x}) \to \phi(\mathbf{x}) \text{ as } h \to 0 \; \forall \phi(\mathbf{x}) \in L^2(\mathcal{D}), \qquad (23)$$

then (21) follows from the compactness of $\mathcal{G}_\Gamma$ on $L^2(\mathcal{D})$ as stated in Proposition 4 already. Consequently, $\|\phi(\mathbf{x}) - \phi_h(\mathbf{x})\|_{L^2(\mathcal{D})} \to 0$ as $h \to 0$, demonstrating convergence of Galerkin solutions. However, the speed of convergence depends on $\mathcal{S}_h$ and the smoothness of unknown solution $\phi(\mathbf{x})$.

### 4.2. Isogeometric Approximation

The Galerkin discretization discussed in the preceding subsection is general because the finite-dimensional subspaces of $L^2(\mathcal{D})$ have yet to be specified. More often than not, the finite-element subspaces are chosen [7], although the subspaces from mesh-free analysis have also been employed [10]. In this work, the subspaces derived from B-splines and NURBS functions, which are the building blocks of IGA, are proposed.

**Theorem 8.** *Let the Galerkin variational form be as described in (19) or (20). Given the space of NURBS functions $\bar{R}_{\mathbf{i},\mathbf{p}}(\mathbf{x})$, $\mathbf{i} \in \mathcal{I}$, $\mathbf{p} \in \mathbb{N}_0^d$, select*

$$\mathcal{S}_h = \mathcal{V}_h = \text{span}\{\bar{R}_{\mathbf{i},\mathbf{p}}\}_{\mathbf{i} \in \mathcal{I}} \subset L^2(\mathcal{D}), \; h > 0,$$

*as the finite-dimensional subspaces, where $\mathcal{V}_h$ is defined in (7). Then the eigenpairs of the variational form are obtained by solving the linear matrix eigenvalue problem*

$$\mathbf{A}\mathbf{f}_h = \lambda_h \mathbf{B}\mathbf{f}_h, \qquad (24)$$

*yielding an eigenvalue $\lambda_h \in \mathbb{R}_0^+$ and an eigenvector $\mathbf{f}_h \in \mathbb{R}^{n_c}$, where $n_c$ is the number of control points of IGA. Here, $\mathbf{A} \in \mathbb{R}^{n_c \times n_c}$ and $\mathbf{B} \in \mathbb{R}^{n_c \times n_c}$ are system matrices, which have components*

$$A_{\mathbf{ij}} := \int_{\mathcal{D}} \int_{\mathcal{D}} \Gamma(\mathbf{x},\mathbf{x}')\bar{R}_{\mathbf{i},\mathbf{p}}(\mathbf{x})\bar{R}_{\mathbf{j},\mathbf{p}}(\mathbf{x}')d\mathbf{x}d\mathbf{x}', \quad \mathbf{i},\mathbf{j} \in \mathcal{I}, \qquad (25)$$

*and*

$$B_{\mathbf{ij}} := \int_{\mathcal{D}} \bar{R}_{\mathbf{i},\mathbf{p}}(\mathbf{x})\bar{R}_{\mathbf{j},\mathbf{p}}(\mathbf{x})d\mathbf{x}, \quad \mathbf{i},\mathbf{j} \in \mathcal{I}, \qquad (26)$$

*with $\mathbf{i},\mathbf{j}$ representing any two control points of IGA. Henceforth, the corresponding eigenfunction is obtained as*

$$\phi_h(\mathbf{x}) = \sum_{\mathbf{j} \in \mathcal{I}} f_{h,\mathbf{j}} \bar{R}_{\mathbf{j},\mathbf{p}}(\mathbf{x}), \qquad (27)$$

*where $f_{h,\mathbf{j}}$ is the $\mathbf{j}$th component of $\mathbf{f}_h \in \mathbb{R}^{n_c}$.*

*Proof.* Since $\psi(\mathbf{x}) \in \mathcal{V}_h$, expand the function

$$\psi(\mathbf{x}) = \sum_{\mathbf{i} \in \mathcal{I}} a_{\mathbf{i}} \bar{R}_{\mathbf{i},\mathbf{p}}(\mathbf{x}) \qquad (28)$$

with respect to the basis of $\mathcal{V}_h \subset L^2(\mathcal{D})$, where $a_{\mathbf{i}}$, $\mathbf{i} \in \mathcal{I}$, are the associated coefficients. Applying (27) and (28) into (20) and interchanging the integral and summation operators gives

$$\sum_{\mathbf{i} \in \mathcal{I}} \sum_{\mathbf{j} \in \mathcal{I}} a_{\mathbf{i}} f_{h,\mathbf{j}} \int_{\mathcal{D}} \int_{\mathcal{D}} \Gamma(\mathbf{x},\mathbf{x}')\bar{R}_{\mathbf{i},\mathbf{p}}(\mathbf{x})\bar{R}_{\mathbf{j},\mathbf{p}}(\mathbf{x}')d\mathbf{x}d\mathbf{x}' = \lambda_h \sum_{\mathbf{i} \in \mathcal{I}} \sum_{\mathbf{j} \in \mathcal{I}} a_{\mathbf{i}} f_{h,\mathbf{j}} \int_{\mathcal{D}} \bar{R}_{\mathbf{i},\mathbf{p}}(\mathbf{x})\bar{R}_{\mathbf{i},\mathbf{p}}(\mathbf{x})d\mathbf{x}. \qquad (29)$$

From the definitions of the system matrices $\mathbf{A}$ and $\mathbf{B}$ as in (25) and (26), (29) becomes

$$\sum_{\mathbf{i}\in\mathcal{I}} a_{\mathbf{i}} \left( \sum_{\mathbf{j}\in\mathcal{I}} A_{\mathbf{ij}} f_{h,\mathbf{j}} - \lambda_h \sum_{\mathbf{j}\in\mathcal{I}} B_{\mathbf{ij}} f_{h,\mathbf{j}} \right) = 0. \tag{30}$$

Since $\psi(\mathbf{x})$ is an arbitrary member of $\mathcal{V}_h$, the constants $a_{\mathbf{i}}$, $\mathbf{i} \in \mathcal{I}$, are also arbitrary. Therefore, the parenthetical term of (30) must vanish for all $\mathbf{i} \in \mathcal{I}$, resulting in (24). □

An alternative proof of Theorem 8 involves the following steps. Define the residual error

$$e_{\mathcal{I}} := \sum_{\mathbf{j}\in\mathcal{I}} f_{h,\mathbf{j}} \left( \int_{\mathcal{D}} \Gamma(\mathbf{x},\mathbf{x}') \bar{R}_{\mathbf{j},\mathbf{p}}(\mathbf{x}') d\mathbf{x}' - \lambda_h \bar{R}_{\mathbf{j},\mathbf{p}}(\mathbf{x}) \right)$$

committed by an isogeometric approximation using the subspace $\mathcal{V}_h$ and associated control points of $\mathcal{I}$. Seek the coefficients $f_{h,\mathbf{j}}$, $\mathbf{j} \in \mathcal{I}$, by recognizing the error to be orthogonal to the subspace $\mathcal{V}_h$. Indeed, setting $(e_{\mathcal{I}}, \bar{R}_{\mathbf{i},\mathbf{p}}(\mathbf{x}))_{L^2(\mathcal{D})} = 0$ for all $\mathbf{i} \in \mathcal{I}$ produces (24).

Does the Galerkin isogeometric method for solving the integral eigenvalue problem converge? The question can be readily answered using Theorem 1, which demonstrates that the sequence of isogeometric approximations using NURBS functions converges for any square-integrable function on $\mathcal{D}$. In this case, the projection $P_h : L^2(\mathcal{D}) \to \mathcal{V}_h$ satisfies the condition in (23). Therefore, the eigensolutions from IGA addressed in this work should converge as discussed in Subsection 4.1. Furthermore, a NURBS subspace is usually endowed with smooth functions, depending on the polynomial order of the underlying B-splines and the multiplicity of knots. In consequence, the results of IGA are expected to be smoother than those of FEM, especially when concerned with continuity across element boundaries. The smoothness of IGA solutions will be further examined in the following section.

4.3. Properties of System Matrices

The solution of the matrix eigenvalue value problem described in Theorem 8 depends on the properties of the system matrices $\mathbf{A}$ and $\mathbf{B}$. The following proposition demonstrates that both matrices are symmetric and are either positive-semidefinite or positive-definite.

**Proposition 9.** Let the system matrices $\mathbf{A}$ and $\mathbf{B}$ be as defined in (25) and (26), respectively. Then $\mathbf{A}$ is symmetric and positive-semidefinite, and $\mathbf{B}$ is symmetric and positive-definite.

*Proof.* From their definitions and the fact that the covariance function $\Gamma(\mathbf{x}, \mathbf{x}')$ is symmetric with respect to its arguments $\mathbf{x}$ and $\mathbf{x}'$, $\mathbf{A}$ and $\mathbf{B}$ are both symmetric matrices. To prove positive-semidefiniteness of $\mathbf{A}$, let $\mathbf{0} \neq \mathbf{c} \in \mathbb{R}^{n_c}$ be a column vector of arbitrary, but non-*zero*, constants $c_{\mathbf{i}}$, $\mathbf{i} \in \mathcal{I}$. Applying (11) from Mercer's theorem and interchanging summation and integrals operators, including the convergent infinite sum,

$$\begin{aligned}
\mathbf{c}^T \mathbf{A} \mathbf{c} &= \sum_{\mathbf{i}\in\mathcal{I}} \sum_{\mathbf{j}\in\mathcal{I}} c_{\mathbf{i}} c_{\mathbf{j}} A_{\mathbf{ij}} \\
&= \sum_{\mathbf{i}\in\mathcal{I}} \sum_{\mathbf{j}\in\mathcal{I}} c_{\mathbf{i}} c_{\mathbf{j}} \int_{\mathcal{D}} \int_{\mathcal{D}} \Gamma(\mathbf{x},\mathbf{x}') \bar{R}_{\mathbf{i},\mathbf{p}}(\mathbf{x}) \bar{R}_{\mathbf{j},\mathbf{p}}(\mathbf{x}') d\mathbf{x} d\mathbf{x}' \\
&= \int_{\mathcal{D}} \int_{\mathcal{D}} \sum_{k=1}^{\infty} \lambda_k \phi_k(\mathbf{x}) \phi_k(\mathbf{x}') \sum_{\mathbf{i}\in\mathcal{I}} c_{\mathbf{i}} \bar{R}_{\mathbf{i},\mathbf{p}}(\mathbf{x}) \sum_{\mathbf{j}\in\mathcal{I}} c_{\mathbf{j}} \bar{R}_{\mathbf{j},\mathbf{p}}(\mathbf{x}') d\mathbf{x} d\mathbf{x}' \\
&= \sum_{k=1}^{\infty} \lambda_k \int_{\mathcal{D}} \int_{\mathcal{D}} \sum_{\mathbf{i}\in\mathcal{I}} c_{\mathbf{i}} \bar{R}_{\mathbf{i},\mathbf{p}}(\mathbf{x}) \phi_k(\mathbf{x}) d\mathbf{x} \sum_{\mathbf{j}\in\mathcal{I}} c_{\mathbf{j}} \bar{R}_{\mathbf{j},\mathbf{p}}(\mathbf{x}') \phi_k(\mathbf{x}') d\mathbf{x}' \\
&= \sum_{k=1}^{\infty} \lambda_k \left( \int_{\mathcal{D}} \sum_{\mathbf{i}\in\mathcal{I}} c_{\mathbf{i}} \bar{R}_{\mathbf{i},\mathbf{p}}(\mathbf{x}) \phi_k(\mathbf{x}) d\mathbf{x} \right)^2 \\
&\geq 0,
\end{aligned}$$

as $\lambda_k$ is a non-negative real number while the square of the integral in the third line is non-negative. Therefore, $\mathbf{A}$ is positive-semidefinite.

Finally, for the matrix $\mathbf{B}$, swapping, again, the summation and integral operators,

$$\begin{aligned}
\mathbf{c}^T \mathbf{B} \mathbf{c} &= \sum_{\mathbf{i} \in \mathcal{I}} \sum_{\mathbf{j} \in \mathcal{I}} c_{\mathbf{i}} c_{\mathbf{j}} B_{\mathbf{ij}} \\
&= \sum_{\mathbf{i} \in \mathcal{I}} \sum_{\mathbf{j} \in \mathcal{I}} c_{\mathbf{i}} c_{\mathbf{j}} \int_{\mathcal{D}} \bar{R}_{\mathbf{i},\mathbf{p}}(\mathbf{x}) \bar{R}_{\mathbf{j},\mathbf{p}}(\mathbf{x}) d\mathbf{x} \\
&= \int_{\mathcal{D}} \sum_{\mathbf{i} \in \mathcal{I}} c_{\mathbf{i}} \bar{R}_{\mathbf{i},\mathbf{p}}(\mathbf{x}) \sum_{\mathbf{j} \in \mathcal{I}} c_{\mathbf{j}} \bar{R}_{\mathbf{j},\mathbf{p}}(\mathbf{x}) d\mathbf{x} \\
&= \left( \sum_{\mathbf{i} \in \mathcal{I}} c_{\mathbf{i}} \bar{R}_{\mathbf{i},\mathbf{p}}(\mathbf{x}), \sum_{\mathbf{j} \in \mathcal{I}} c_{\mathbf{j}} \bar{R}_{\mathbf{j},\mathbf{p}}(\mathbf{x}) \right)_{L^2(\mathcal{D})} \\
&= \left\| \sum_{\mathbf{i} \in \mathcal{I}} c_{\mathbf{i}} \bar{R}_{\mathbf{i},\mathbf{p}}(\mathbf{x}) \right\|^2_{L^2(\mathcal{D})} \\
&> 0,
\end{aligned}$$

as the norm is positive for any $\mathbf{c} \neq \mathbf{0}$, thereby proving its positive-definiteness. □

Proposition 9 is applicable for any polynomial order $p$ of NURBS. In fact, the symmetry and positive-(semi)definiteness are valid when the system matrices are defined in conjunction with other subspaces rooted in finite-element and mesh-free analyses. However, the definitions of respective system matrices depend on the basis functions of individual subspaces.

Given the results of Proposition 9, the following properties of eigensolutions are fulfilled: (1) all eigenvalues are real and non-negative; (2) the eigenvectors corresponding to distinct eigenvalues are mutually orthogonal; (3) for repeated eigenvalues, $m \in \mathbb{N}$ mutually orthogonal eigenvectors can be found for each eigenvalue of multiplicity $m$; and (4) the eigenvectors constitute a basis of $\mathbb{R}^{n_c}$. These properties are consistent with those endowed to the Hilbert-Schmidt operator $\mathcal{G}_\Gamma : L^2(\mathcal{D}) \to L^2(\mathcal{D})$, as explained by Proposition 5.

The size of both system matrices is $n_c \times n_c$, where $n_c := |\mathcal{I}|$ is the number of control points defined in Section 2. Obviously, the largest number of eigensolutions of (24) is limited to $n_c$. This is not an issue for most IGA meshes, as $n_c$ is typically larger than $N$, the number of eigensolutions retained in the KL approximation. On the contrary, if $n_c < N$, then all needed eigensolutions cannot be obtained. This will be further clarified in Section 5 where numerical results are discussed.

### 4.4. Construction of System Matrices

The assembly of the system matrices $\mathbf{A}$ and $\mathbf{B}$ requires domain integrations in the physical space. In general, these integrals cannot be determined exactly. Therefore, the matrices must be estimated by numerical integration. However, the use of NURBS functions in isogeometric analysis introduces the parametric domain as explained in Section 2. This slightly complicates the matter because, for numerical integration, an additional domain $[-1, +1]^d$ is needed. The latter domain is commonly referred to as the parent element in the isogeometric literature [12].

Consider an arbitrary element $Q \in \mathcal{Q}_h$ in the parametric domain $\hat{\mathcal{D}} = [0,1]^d$. Each such element can be viewed as the image of the parent element $[-1, +1]^d$ defined by an affine mapping $\boldsymbol{\xi} : [-1, +1]^d \to Q$ or, equivalently, by $\boldsymbol{\xi}(\boldsymbol{\eta})$, where $\boldsymbol{\eta}$ is the coordinate of the parent element. Similarly, there is a corresponding element $K \in \mathcal{K}_h$ in the physical domain $\mathcal{D} \subset \mathbb{R}^d$, which is the image of that very element $Q$ in the parametric domain. Recall that $\mathbf{x} : \hat{\mathcal{D}} \to \mathcal{D}$, that is, $\mathbf{x}(\boldsymbol{\xi})$, is the mapping between the parametric and physical domains. The same mapping is used for these two corresponding elements. To integrate a function on an element of the physical domain, a pullback of the physical element to the parent element is required. This is accomplished using the composition of the inverses, that is, $\mathbf{x}^{-1} : K \to Q$ and $\boldsymbol{\xi}^{-1} : Q \to [-1, +1]^d$, of the two aforementioned mappings.

Let the Jacobians of the mappings $\boldsymbol{\xi}(\boldsymbol{\eta})$ and $\mathbf{x}(\boldsymbol{\xi})$ be defined as

$$\mathbf{J}_{\boldsymbol{\eta}} := \left[\frac{\partial \boldsymbol{\xi}}{\partial \boldsymbol{\eta}}\right] = \begin{bmatrix} \frac{\partial \xi_1}{\partial \eta_1} & \cdots & \frac{\partial \xi_1}{\partial \eta_d} \\ \vdots & \ddots & \vdots \\ \frac{\partial \xi_d}{\partial \eta_1} & \cdots & \frac{\partial \xi_d}{\partial \eta_d} \end{bmatrix}$$

and

$$\mathbf{J}_{\boldsymbol{\xi}} := \left[\frac{\partial \mathbf{x}}{\partial \boldsymbol{\xi}}\right] = \begin{bmatrix} \frac{\partial x_1}{\partial \xi_1} & \cdots & \frac{\partial x_1}{\partial \xi_d} \\ \vdots & \ddots & \vdots \\ \frac{\partial x_d}{\partial \xi_1} & \cdots & \frac{\partial x_d}{\partial \xi_d} \end{bmatrix},$$

respectively. As the mapping is affine, the calculation of the partial derivatives $\partial \xi_i / \partial \eta_j$, $i, j = 1, \ldots, d$, is straightforward. However, to determine the partial derivatives $\partial x_i / \partial \xi_j$, $i, j = 1, \ldots, d$, the derivatives of NURBS and B-Spline functions are involved. Due to brevity, explicit details of the derivatives of NURBS and B-spline functions are not reported here as they are available elsewhere [12].

Given the mappings and their respective Jacobians, the components of the system matrices are then evaluated by summing contributions from all element-level integrations on the parent element, that is,

$$\begin{aligned}
A_{\mathbf{ij}} &:= \int_{\mathcal{D}} \int_{\mathcal{D}} \Gamma(\mathbf{x}, \mathbf{x}') \bar{R}_{\mathbf{i},\mathbf{p}}(\mathbf{x}) \bar{R}_{\mathbf{j},\mathbf{p}}(\mathbf{x}') d\mathbf{x} d\mathbf{x}' \\
&= \sum_{K \in \mathcal{K}_h} \sum_{K' \in \mathcal{K}_h} \int_K \int_{K'} \Gamma(\mathbf{x}, \mathbf{x}') \bar{R}_{\mathbf{i},\mathbf{p}}(\mathbf{x}) \bar{R}_{\mathbf{j},\mathbf{p}}(\mathbf{x}') d\mathbf{x} d\mathbf{x}' \\
&= \sum_{Q \in \mathcal{Q}_h} \sum_{Q' \in \mathcal{Q}_h} \int_Q \int_{Q'} \Gamma(\mathbf{x}(\boldsymbol{\xi}), \mathbf{x}'(\boldsymbol{\xi}')) \bar{R}_{\mathbf{i},\mathbf{p}}(\mathbf{x}(\boldsymbol{\xi})) \bar{R}_{\mathbf{j},\mathbf{p}}(\mathbf{x}'(\boldsymbol{\xi}')) |\det \mathbf{J}_{\boldsymbol{\xi}}| d\boldsymbol{\xi} d\boldsymbol{\xi}' \\
&= \sum \sum \int_{[-1,+1]^d} \int_{[-1,+1]^d} \Gamma(\mathbf{x}(\boldsymbol{\xi}(\boldsymbol{\eta})), \mathbf{x}'(\boldsymbol{\xi}'(\boldsymbol{\eta}'))) \bar{R}_{\mathbf{i},\mathbf{p}}(\mathbf{x}(\boldsymbol{\xi}(\boldsymbol{\eta}))) \bar{R}_{\mathbf{j},\mathbf{p}}(\mathbf{x}'(\boldsymbol{\xi}'(\boldsymbol{\eta}'))) |\det \mathbf{J}_{\boldsymbol{\xi}}| |\det \mathbf{J}_{\boldsymbol{\eta}}| d\boldsymbol{\eta} d\boldsymbol{\eta}'
\end{aligned} \quad (31)$$

and

$$\begin{aligned}
B_{\mathbf{ij}} &:= \int_{\mathcal{D}} \bar{R}_{\mathbf{i},\mathbf{p}}(\mathbf{x}) \bar{R}_{\mathbf{j},\mathbf{p}}(\mathbf{x}) d\mathbf{x} \\
&= \sum_{K \in \mathcal{K}_h} \int_K \bar{R}_{\mathbf{i},\mathbf{p}}(\mathbf{x}) \bar{R}_{\mathbf{j},\mathbf{p}}(\mathbf{x}) d\mathbf{x} \\
&= \sum_{Q \in \mathcal{Q}_h} \int_Q \bar{R}_{\mathbf{i},\mathbf{p}}(\mathbf{x}(\boldsymbol{\xi})) \bar{R}_{\mathbf{j},\mathbf{p}}(\mathbf{x}(\boldsymbol{\xi})) |\det \mathbf{J}_{\boldsymbol{\xi}}| d\boldsymbol{\xi} \\
&= \sum \int_{[-1,+1]^d} \bar{R}_{\mathbf{i},\mathbf{p}}(\mathbf{x}(\boldsymbol{\xi}(\boldsymbol{\eta}))) \bar{R}_{\mathbf{j},\mathbf{p}}(\mathbf{x}(\boldsymbol{\xi}(\boldsymbol{\eta})))|\det \mathbf{J}_{\boldsymbol{\xi}}| |\det \mathbf{J}_{\boldsymbol{\eta}}| d\boldsymbol{\eta}.
\end{aligned} \quad (32)$$

Here, the summations in the last lines of (31) and (32) are over all $n_e := |\mathcal{Q}_h| = |\mathcal{K}_h|$ elements of IGA. The final integrals in (31) and (32) are estimated by a suitable numerical integration scheme, such as the Gauss-Legendre quadrature. Even though the NURBS functions are not necessarily polynomials, the Gauss quadrature is still effective [12]. In this case, the same quadrature rule employed for a $p$th-order polynomial can be used for a NURBS function built from an underlying $p$th-order B-spline. Having said so, a Gauss quadrature is not an optimal choice for IGA. That is why current research is focused on finding optional or near-optimal numerical integration techniques to tackle NURBS functions [24].

It is important to underscore that the construction of the system matrices $\mathbf{A}$ and $\mathbf{B}$ mandates, respectively, $2d$- and $d$-dimensional domain integrations – a fundamental imposition of the Galerkin discretization. While forming $\mathbf{B}$ requires an effort similar to that of assembling the mass matrix in solid mechanics, building $\mathbf{A}$ is computationally daunting as it is $d$-order more expensive than forming $\mathbf{B}$. For instance, on three-dimensional ($d = 3$) domains, $\mathbf{A}$ requires a six-fold integration as opposed to a three-fold integration needed by **B**. To speed up the assembly, a few researchers [25, 26] have suggested replacing **A** with a hierarchical matrix [27], but this also adds a new layer of approximation. Nonetheless, the implementation of an industrial-scale matrix eigenvalue problem in high dimension is a formidable task. Such a computational challenge is not unique to IGA; it persists in finite-element and mesh-free analyses as well.

*4.5. Eigenvalue Solvers*

While the symmetry and positive-(semi)definiteness of the system matrices facilitate seeking only real-valued eigensolutions, the eigenvalue problem presents a few computational challenges of its own. First, depending on how fast or how slow the covariance function decays, the number of eigenmodes retained in a truncated KL expansion can be large or small. If, indeed, the required number of eigenmodes is very large, especially for a problem with complex domain encompassing a large number of elements or control points, the computational effort in deriving the concomitant eigensolutions can also be very large. Second, in many cases, **B** is a relatively sparse and diagonalizable matrix. In contrast, **A** is a dense matrix, making its generation and storage expensive. That is why selecting efficient eigenvalue solvers, such as Krylov subspace iteration in conjunction with the Arnoldi method [28], the Lanczos-based thick-restart method [29], and the fast multipole method [30], continues to be studied by some researchers [22, 25, 26].

## 5. Numerical Examples

Three examples, each with a progressively higher dimension, are presented to illustrate the proposed Galerkin isogeometric method for solving the integral eigenvalue problem associated with the KL expansion. The random fields have *zero* means and are homogeneous and isotropic in all examples. The polynomial order $p$ of B-splines, leading to NURBS, and the order of Gauss-Legendre quadrature to estimate the system matrices vary from example to example. The isogeometric analysis and subsequent matrix eigenvalue calculations were performed using MATLAB (Version 2016a) [31]. The eigenvalue calculations were checked using multiple algorithms and methods, such as Cholesky factorization, the QZ algorithm, and the Lanczos method, all of which are available in MATLAB.

*5.1. Example 1*

Consider a one-dimensional random field $\alpha(x,\cdot)$ with covariance function $\Gamma(x,x') = \mathbb{E}[\alpha(x,\cdot)\alpha(x',\cdot)]$ defined on $\mathcal{D} = [0,1] \subset \mathbb{R}$. Three types of covariance functions, described by exponential, Gaussian, and Bessel functions, were selected. Mathematically,

$$\Gamma(x,x') = \begin{cases} \sigma^2 \exp\left(-\dfrac{|x-x'|}{bL}\right), & \text{Type 1 (exponential)}, \\ \sigma^2 \exp\left(-\dfrac{|x-x'|^2}{(bL)^2}\right), & \text{Type 2 (Gaussian)}, \\ \sigma^2 \dfrac{|x-x'|}{bL} K_1\left(\dfrac{|x-x'|}{bL}\right), & \text{Type 3 (Bessel)}, \end{cases}$$

where $K_1$ is the modified second-kind Bessel function of order one, $\sigma^2 = 1$, $L = 1$, and $b = 1$ or $0.1$ as the correlation length parameter. These covariance functions, depicted in Figure 2, are commonly used in engineering and applied sciences.

Three polynomial orders, $p=1$, $p=2$, and $p=3$, representing, respectively, linear, quadratic, and cubic elements, were employed. The knot vectors for the coarsest one-element IGA meshes are as follows: (1) $\Xi = (0,0,1,1)$ for linear elements; (2) $\Xi = (0,0,0,1,1,1)$ for quadratic elements; and (3) $\Xi = (0,0,0,0,1,1,1,1)$ for cubic elements. The associated control points and weights are given in Table A.1 of Appendix A. As the weights are all equal to one, the NURBS functions are the same as B-splines. By adding new knots and control points, the one-element mesh for each polynomial order was $h$-refined globally, resulting in a series of progressively finer meshes with the number of elements increasing in size. Two mesh sizes were examined: (1) a coarse mesh with 16 elements, that is, $n_e=16$; and (2) a fine mesh with 256 elements, that is, $n_e=256$.

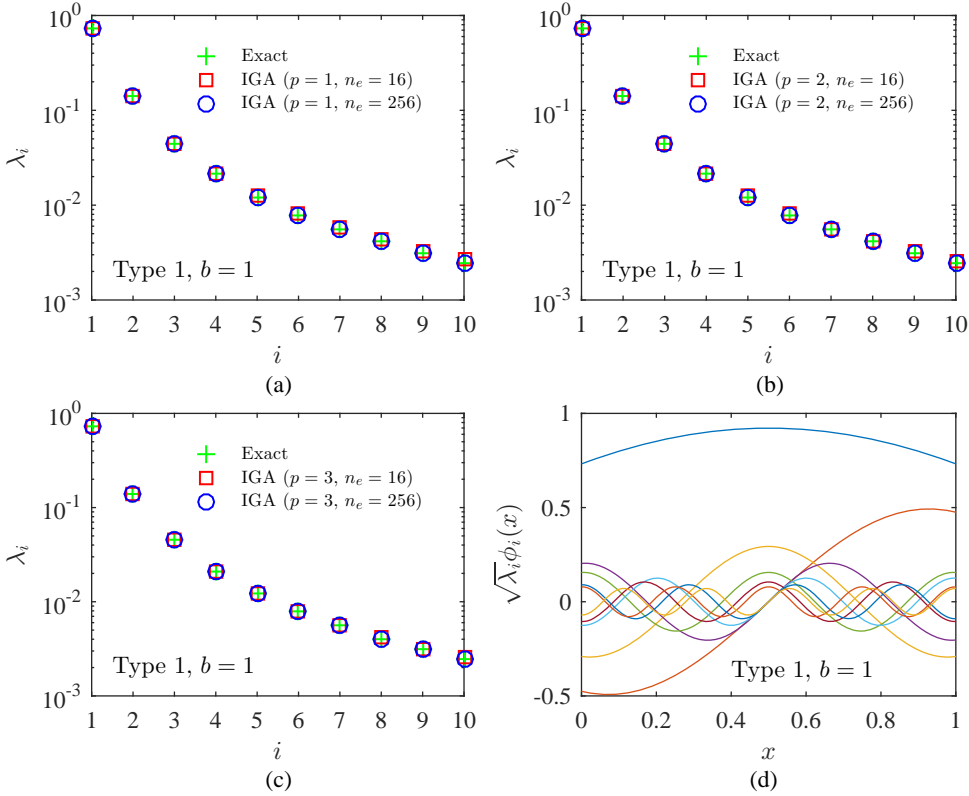

Figure 3: Ten largest exact and IGA-derived eigensolutions for the exponential covariance function with $b = 1$ in Example 1; (a) exact eigenvalues and linear IGA solutions; (b) exact eigenvalues and quadratic IGA solutions; (c) exact eigenvalues and cubic IGA solutions; (d) exact eigenfunctions.

The number of control points varies accordingly and fittingly with the order of elements. A $(p + 1)$-point Gauss-Legendre quadrature was employed to estimate the system matrices and produce the IGA results of Figures 3 through 5.

Figures 3(a) through 3(c) present the scatter plots of the first ten largest eigenvalues for the exponential (Type 1) covariance function with $b = 1$. The exact eigenvalues, which exist for this covariance function [5, 7], and approximate eigenvalues obtained using the proposed isogeometric method are displayed. The eigenvalue decays with the mode number as expected. For the IGA results, the figures are arranged according to linear [Figure 3(a)], quadratic [Figure 3(b)], and cubic [Figure 3(c)] elements, each containing the results of the coarse mesh (16 elements) and the fine mesh (256 elements), as defined previously. The quality of agreement between the exact eigenvalues and their respective IGA approximations is excellent even for the coarse mesh. Any distinction between the exact and IGA solutions is impalpable to the naked eye, especially when examining the results of higher-order IGA.

The same observation holds when comparing the exact (scaled) eigenfunctions, also available for the exponential covariance function [5, 7], in Figure 3(d) and their IGA derived approximations in Figures 4(a) through 4(d) when $b = 1$. Again, the first ten eigenfunctions are shown. The eigenfunctions obtained from IGA are organized according to the use of (1) linear elements with the coarse [Figure 4(a)] and fine [Figure 4(b)] meshes; (2) quadratic elements with the coarse mesh [Figure 4(c)]; and (3) cubic elements with the coarse mesh [Figure 4(d)]. In Figure 4(a), the IGA-derived eigenfunctions show a lack of smoothness, but this is expected for linear elements with the coarse mesh. The smoothness returns, and the results get better for the fine mesh, as confirmed in Figure 4(b). Having said this, the eigenfunctions are still $C^0$-continuous across element boundaries regardless of the mesh size because of linear basis functions. Indeed, the same result should be expected from traditional FEM employing piecewise-linear basis functions. In contrast, when quadratic or cubic elements are used, the eigenfunctions from IGA even for the coarse mesh, shown in Figure 4(c) or 4(d), are already smooth and similar to, if not better than, the results of linear elements with the fine mesh. In addition, these higher-order IGA solutions are endowed with $C^1$- or $C^2$- continuity across element boundaries when compared with low- or high-order FEM with $C^0$-continuity.

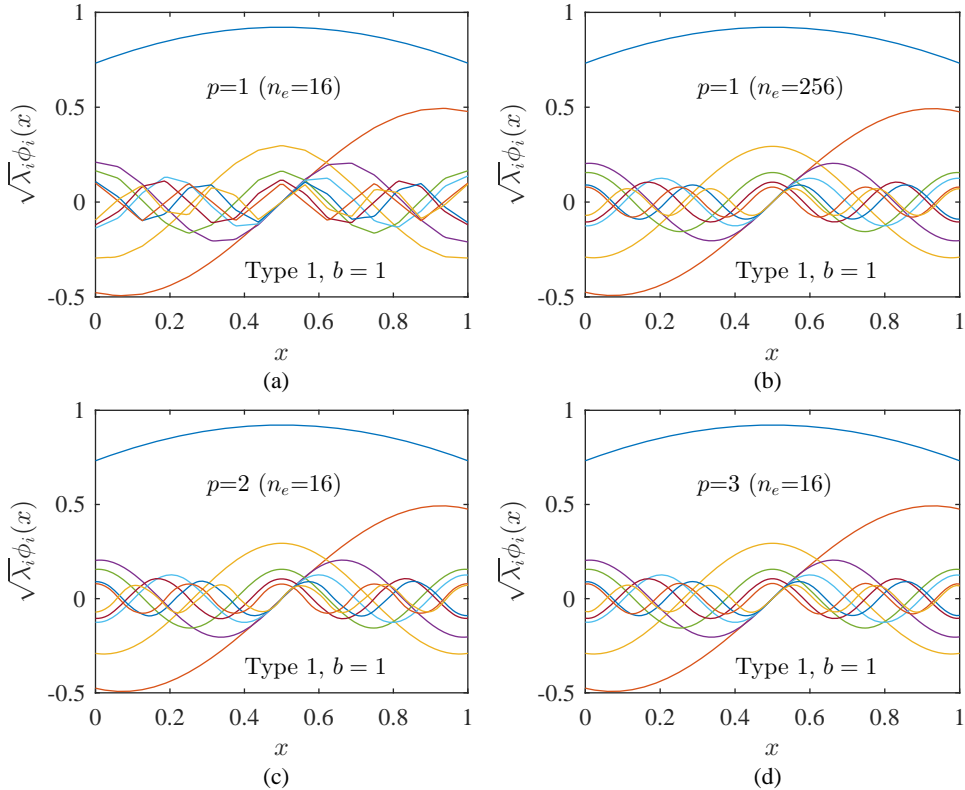

Figure 4: Ten largest IGA-derived eigenfunctions for the exponential covariance function with $b = 1$ in Example 1; (a) linear IGA with 16 elements; (b) linear IGA with 256 elements; (c) quadratic IGA with 16 elements; (d) cubic IGA with 16 elements.

Figures 5(a) and 5(b) contrast the reference and IGA solutions of the first ten eigenvalues for the Gaussian (Type 2) and Bessel (Type 3) covariance functions, respectively, both employing the correlation length parameter $b = 1$. Since no exact solutions exist for these two covariance functions, respective IGA solutions obtained for a further refined mesh of 512 elements were used as reference solutions. In each figure, the top, middle, and bottom sub-figures contain the results of linear, quadratic, and cubic elements, respectively, from both the coarse mesh and the fine mesh. The IGA-derived eigenvalues converge to the reference solutions for both types of covariance functions. No comparisons between reference solutions and IGA results of eigenfunctions are reported here, as they essentially show the same qualitative trend exhibited for the exponential covariance function.

Finally, a limited numerical error analysis of IGA-derived eigenvalues was performed for the Gaussian and exponential covariance functions with $b = 0.1$. More specifically, a relative error, defined as the ratio of (1) the 2-norm of the vector of first ten eigenvalue differences between reference or exact solutions and IGA solutions and (2) the 2-norm of the vector of first ten reference or exact eigenvalues, was studied. To eliminate possible errors due to numerical integration, a higher-order quadrature, that is, a $(p + 5)$-point Gauss-Legendre quadrature was employed for obtaining all IGA solutions, including the reference IGA solution entailing 512 cubic elements for the Gaussian covariance function. The results are delineated in Figures 6(a) and 6(b), explaining how the relative error decays with respect to the element size $h$, which is reciprocal to the number of elements and hence constant for a fixed number of elements. For either covariance function, the error decays with the reduction of the element size, as expected, regardless of the polynomial order $p$. Given a fixed mesh size, the error from a higher-order IGA is consistently lower than that from the low-order IGA. Moreover, the rate of decay in Figure 6(a) is faster for the higher-order IGA for the Gaussian covariance function. This is due to not only IGA, but also the smoothness properties of the Gaussian covariance function. By contrast, no such convergence acceleration is found in Figure 6(b), as the exponential covariance function is a non-differentiable one. Therefore, the regularity of the covariance function plays an important role on the performance of any numerical method, including IGA. A theoretical error analysis of IGA-derived eigensolutions supporting numerical trends merits further study.

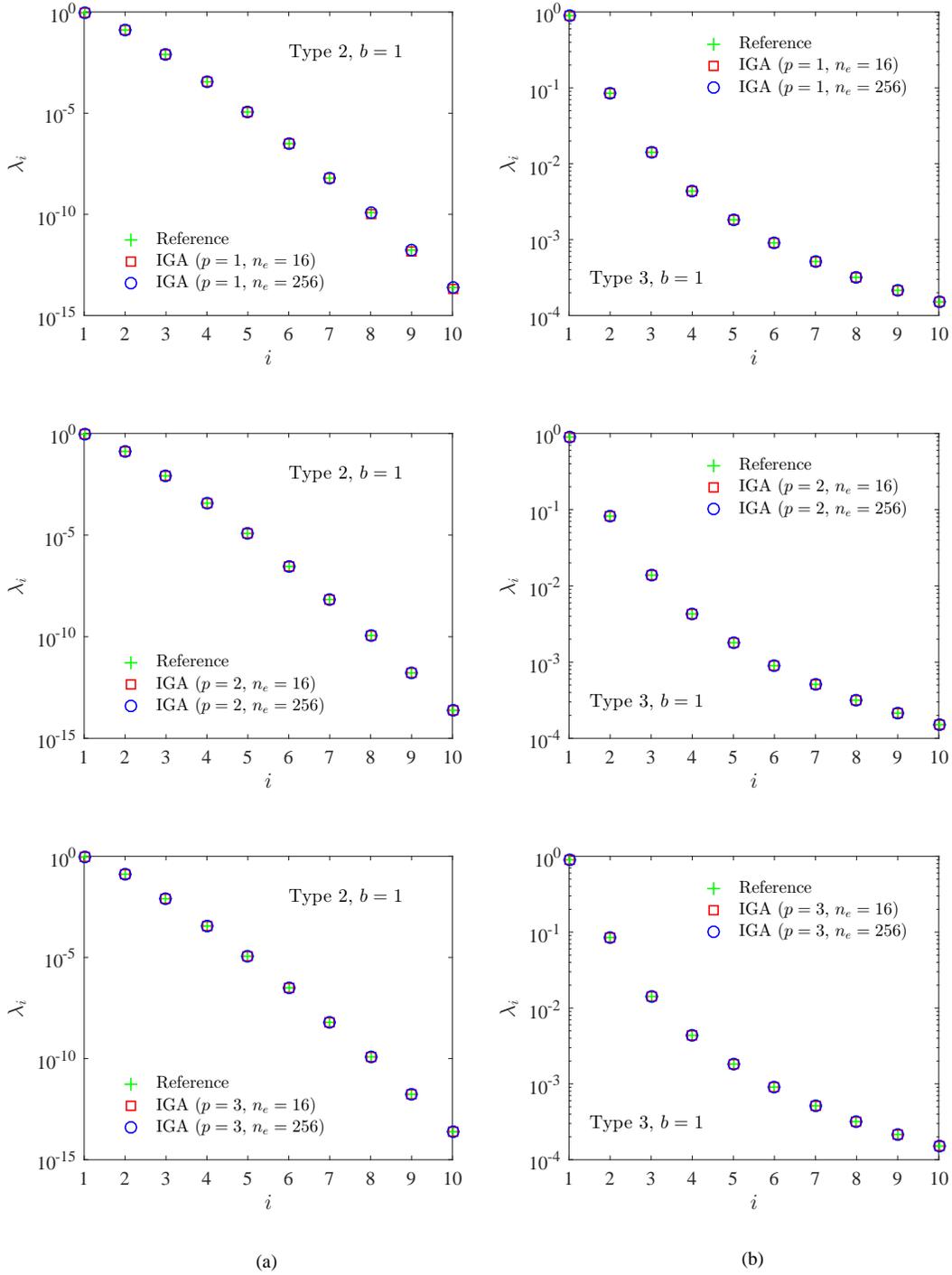

Figure 5: Ten largest eigenvalues from the reference solutions and IGA for the two other covariance functions with $b = 1$ in Example 1; (a) Gaussian covariance function; (b) Bessel covariance function.

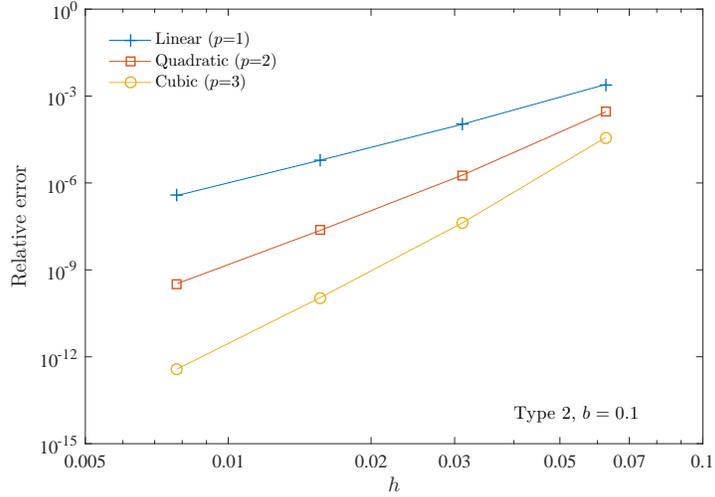

(a)

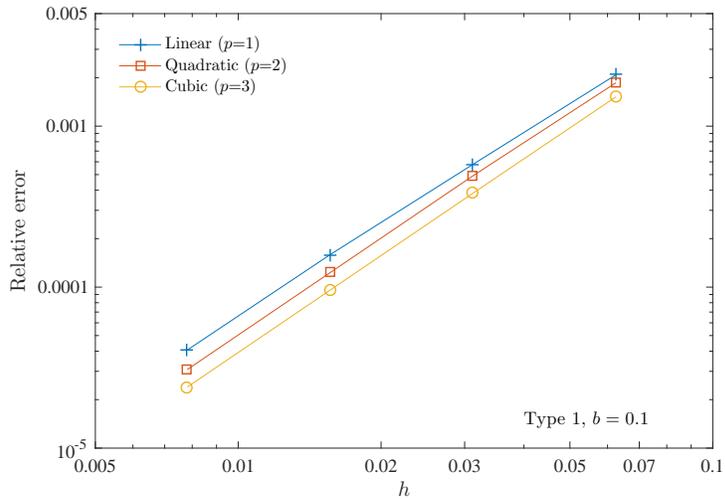

(b)

Figure 6: Relative errors in IGA-derived eigenvalues with $b = 0.1$ in Example 1; (a) Gaussian covariance function; (b) exponential covariance function

## 5.2. Example 2

For the second example, let $\alpha(\mathbf{x}, \cdot)$ be a two-dimensional random field with covariance function

$$\Gamma(\mathbf{x}, \mathbf{x}') = \sigma^2 \exp\left(-\frac{\|\mathbf{x} - \mathbf{x}'\|}{bL}\right), \quad \mathbf{x}, \mathbf{x}' \in \mathcal{D} \subset \mathbb{R}^2, \tag{33}$$

defined on a quarter-annulus of inner radius $R_i = 0.6$ and outer radius $R_o = 1$, as depicted in Figure 7(a). The following covariance parameters were selected: $\sigma^2 = 1$, $L = 1$, and $b = 0.5$.

The initial knot vectors and polynomial orders for the coarsest one-element IGA mesh are as follows: (1) $\Xi_1 = (0, 0, 1, 1)$, $\Xi_2 = (0, 0, 0, 1, 1, 1)$; and (2) $p_1 = 1$, $p_2 = 2$. The initial knot index space and parametric domain are depicted in Figures 7(b) and 7(c), respectively. The associated control points and weights are given in Table A.2 of Appendix A. The weights were chosen in such a way that the resulting NURBS functions could replicate exactly a circular arc. Adopting a global $h$-refinement strategy, the knot index space was successively divided, and new control points and weights were added as needed. Figure 8 displays six meshes, that is, Meshes 1 through 6, with corresponding control points marked in red. The number of elements in these meshes varies from one to 1024. The one-element mesh (Mesh 1), obtained using the initial knot index space, control points, and weights in Table A.2, represents already the exact geometry of the quarter-annulus. In fact, all finer meshes, albeit they have more and more elements, represent the exact geometry of the physical domain. This is in sharp contrast with FEM, where a mesh, especially when it is coarse, will always incur some geometrical errors. A 4-point Gauss-Legendre quadrature was employed in each coordinate direction.

Table 1 presents the first ten largest eigenvalues for the covariance function in (33) calculated by IGA for all six meshes in Figure 8. With the exception of Mesh 1, which has six control points, there are more than ten control points in Meshes 2 through 6. Therefore, all ten eigensolutions were calculated for Meshes 2 through 6, whereas only the first six were calculated for Mesh 1. Clearly, the eigenvalues converge with respect to the number of elements or mesh refinement, as expected. A comparison of the first four eigenfunctions, displayed as contour plots in Figure 9 for Mesh 5 (256 elements) and in Figure 10 for Mesh 6 (1024 elements), tells the same tale with regards to the convergence of eigenfunctions. Note that, unlike in Example 1, there are no analytical eigensolutions for the annular domain and inseparable covariance function in Example 2. That is why multiple meshes were analyzed.

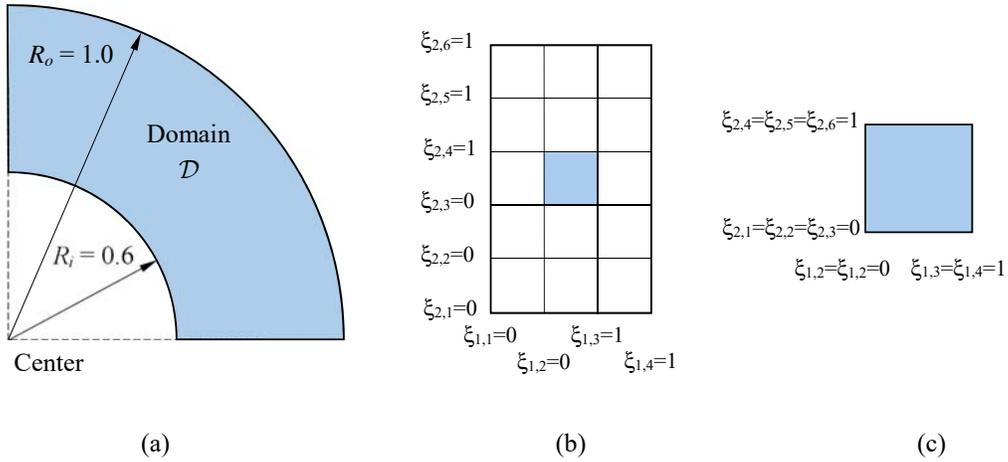

Figure 7: A quarter-annulus in Example 2; (a) physical domain; (b) initial knot index space; (c) parametric domain.

Table 1: Ten largest eigenvalues estimated by IGA using six meshes in Example 2.

| | Eigenvalue | | | | | |
|---|---|---|---|---|---|---|
| | Mesh 1 | Mesh 2 | Mesh 3 | Mesh 4 | Mesh 5 | Mesh 6 |
| Mode | (1 element) | (4 elements) | (16 elements) | (64 elements) | (256 elements) | (1024 elements) |
| 1 | 0.241387791 | 0.234271475 | 0.233414419 | 0.233302062 | 0.233287647 | 0.23328582 |
| 2 | 0.087233705 | 0.086324241 | 0.085394757 | 0.085289015 | 0.085275164 | 0.085273409 |
| 3 | 0.035232598 | 0.032777456 | 0.035032958 | 0.034994641 | 0.034980189 | 0.034978308 |
| 4 | 0.024420835 | 0.020293213 | 0.019896928 | 0.019806854 | 0.019793397 | 0.01979163 |
| 5 | 0.013968298 | 0.016212808 | 0.016475201 | 0.016407932 | 0.016395342 | 0.016393401 |
| 6 | 0.008526274 | 0.012551817 | 0.01211771 | 0.012029456 | 0.012017404 | 0.012015812 |
| 7 | —[a] | 0.006899954 | 0.007744419 | 0.00885116 | 0.008846743 | 0.008844794 |
| 8 | —[a] | 0.004688903 | 0.006573091 | 0.007431735 | 0.007420659 | 0.007419054 |
| 9 | —[a] | 0.004039597 | 0.005035944 | 0.005249849 | 0.005267336 | 0.005265465 |
| 10 | —[a] | 0.003428753 | 0.00449069 | 0.004705134 | 0.004700002 | 0.004698406 |

(a) Not calculable as there are only six control points.

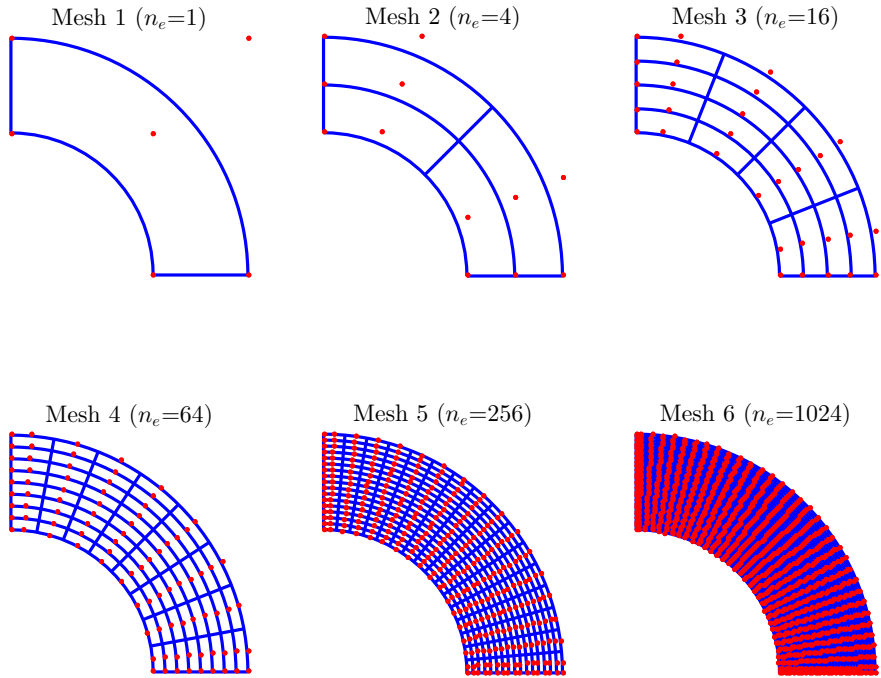

Figure 8: Six IGA meshes obtained by $h$-refinement in Example 2.

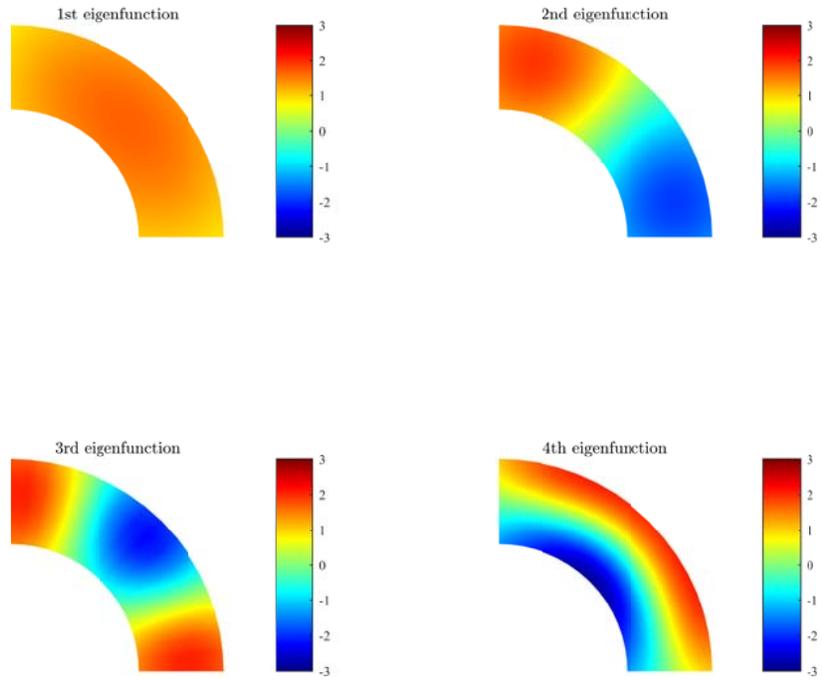

Figure 9: Contour plots of first four IGA-derived eigenfunctions using Mesh 5.

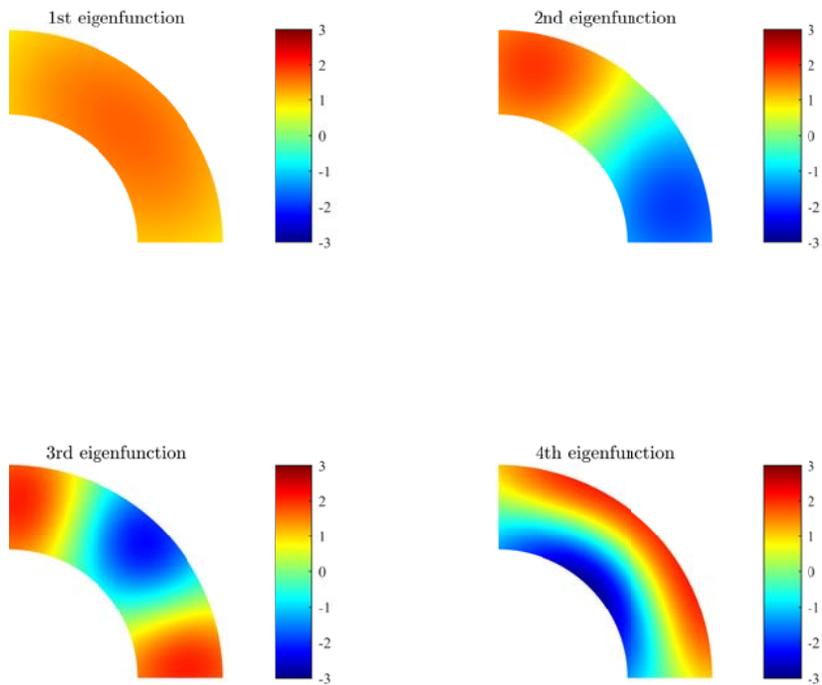

Figure 10: Contour plots of first four IGA-derived eigenfunctions using Mesh 6.

## 5.3. Example 3

The final example entails a three-dimensional random field $\alpha(\mathbf{x}, \cdot)$ with covariance function

$$\Gamma(\mathbf{x}, \mathbf{x}') = \sigma^2 \exp\left(-\frac{\|\mathbf{x} - \mathbf{x}'\|}{bL}\right), \quad \mathbf{x}, \mathbf{x}' \in \mathcal{D} \subset \mathbb{R}^3, \qquad (34)$$

defined on a half-cylinder of inner radius $R_i = 8$, outer radius $R_o = 10$, and length $L_c = 15$, as displayed in Figure 11(a). The covariance parameters were chosen as follows: $\sigma^2 = 1$, $L = 10$ and $b = 0.5$.

Two IGA meshes – one a relatively coarse mesh and the other a relatively fine mesh – were analyzed. The fine mesh was created by doubling the mesh density of the coarse mesh in the circumferential direction, while maintaining the same mesh densities in the radial and length directions. The knot vectors and polynomial orders for both meshes are defined in Table A.3 of Appendix A. However, the numbers of control points and weights, 570 for the coarse mesh and 1050 for the fine mesh, are too many to be listed. The corresponding numbers of elements are 128 and 256, respectively. Both meshes are displayed in Figure 11(b). Again, they represent the exact geometry of the physical domain regardless of the mesh. A 4-point Gauss-Legendre quadrature was employed in each coordinate direction.

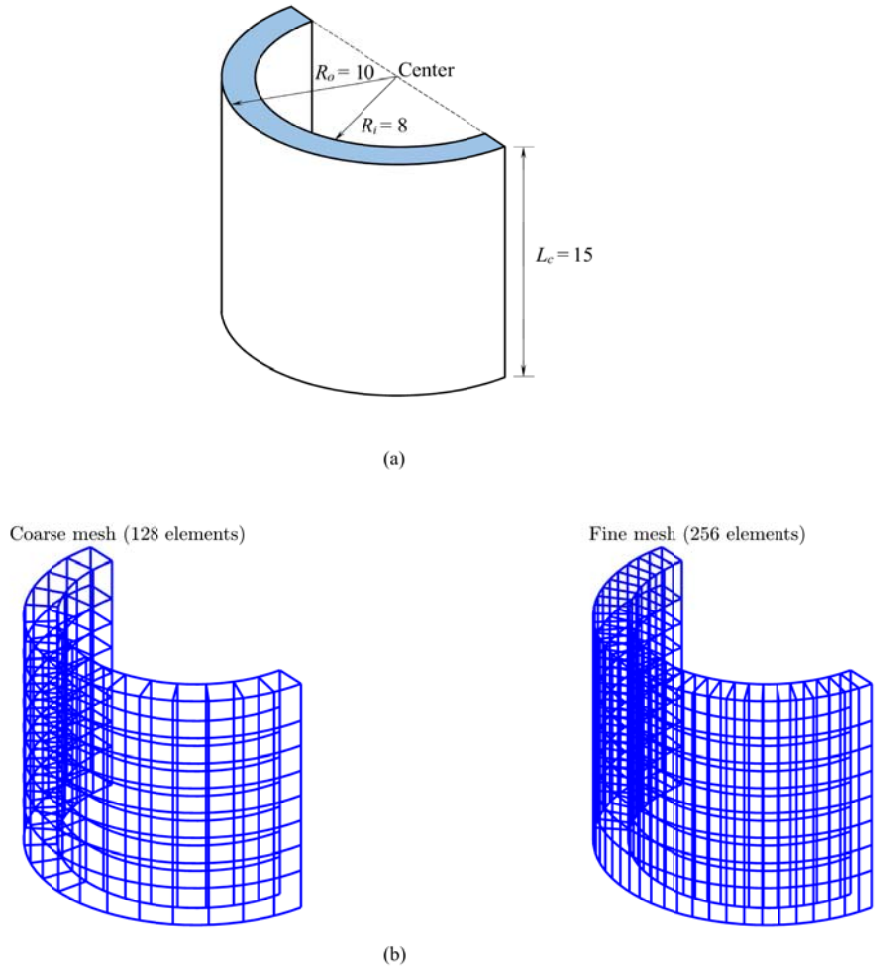

Figure 11: A half-cylinder in Example 3; (a) physical domain; (b) two IGA meshes obtained by circumferential $h$-refinement.

Table 2: Twenty largest eigenvalues estimated by IGA using coarse and fine meshes in Example 3.

| Mode | Eigenvalue | |
| --- | --- | --- |
| | Coarse Mesh (128 elements) | Fine Mesh (256 elements) |
| 1 | 162.8012956 | 162.7993688 |
| 2 | 91.43276902 | 91.43079317 |
| 3 | 57.56984816 | 57.56765769 |
| 4 | 51.09209411 | 51.09029418 |
| 5 | 38.80033743 | 38.79808752 |
| 6 | 27.90574988 | 27.90392169 |
| 7 | 25.05887115 | 25.05681566 |
| 8 | 19.37089182 | 19.36866576 |
| 9 | 16.15870621 | 16.15702078 |
| 10 | 15.79809919 | 15.79601336 |
| 11 | 15.1484597 | 15.1461372 |
| 12 | 11.21556267 | 11.21341529 |
| 13 | 10.1816058 | 10.17964708 |
| 14 | 9.694963725 | 9.693535348 |
| 15 | 8.057339474 | 8.055143387 |
| 16 | 7.579083704 | 7.576970969 |
| 17 | 6.724572177 | 6.722768096 |
| 18 | 6.446973861 | 6.444800351 |
| 19 | 6.178216445 | 6.177288031 |
| 20 | 5.766381043 | 5.764324956 |

Table 2 lists the first twenty largest eigenvalues obtained by IGA using the coarse (128 elements) and fine (256 elements) meshes for the covariance function defined in (34). The respective eigenvalues from both meshes are very close to each other. The same can be said about the eigenfunctions. Due to brevity, however, only the first six eigenfunctions obtained for the fine mesh are portrayed in Figure 12. Again, no analytical solutions exist for this problem, but the relative invariance of eigensolutions from the two meshes provides confidence in the proposed IGA method.

While the computational efforts in the first two examples are relatively low or modest, generating the IGA solution in Example 3, which involves a three-dimensional domain, is computationally demanding. For instance, the total run time for the fine-mesh IGA results in Example 3, obtained using an Intel Core i7-4770, 3.4 Ghz, 16 GB RAM PC, exceeded 24 hours. The root cause for this high computational cost is the double integral over three-dimensional domain in assembling the system matrix $\mathbf{A}$, alluded to in Subsection 4.4. Therefore, developing efficient methods for estimating or constructing $\mathbf{A}$, in the context of the Galerkin framework, is desirable for future work.

## 6. Discussion

While the paper focuses on the Galerkin-based IGA for random field discretization, a brief discussion on the practical implementation of the work and future endeavors is warranted. First, the Galerkin method mandates $2N$-dimensional domain integration for building the system matrix $\mathbf{A}$; this is computationally intensive, as illustrated in the third example, and possibly prohibitive for industrial-scale, three-dimensional applications. Therefore, alternative isogeometric formulations entailing, for instance, collocation methods, which require at most $N$-dimensional domain integrations to construct the system matrices, should be explored to determine their accuracy. Indeed, the collocation method for random field discretization, while exploiting the geometrical flexibility of IGA, is expected to offer a huge computational advantage over the Galerkin method.

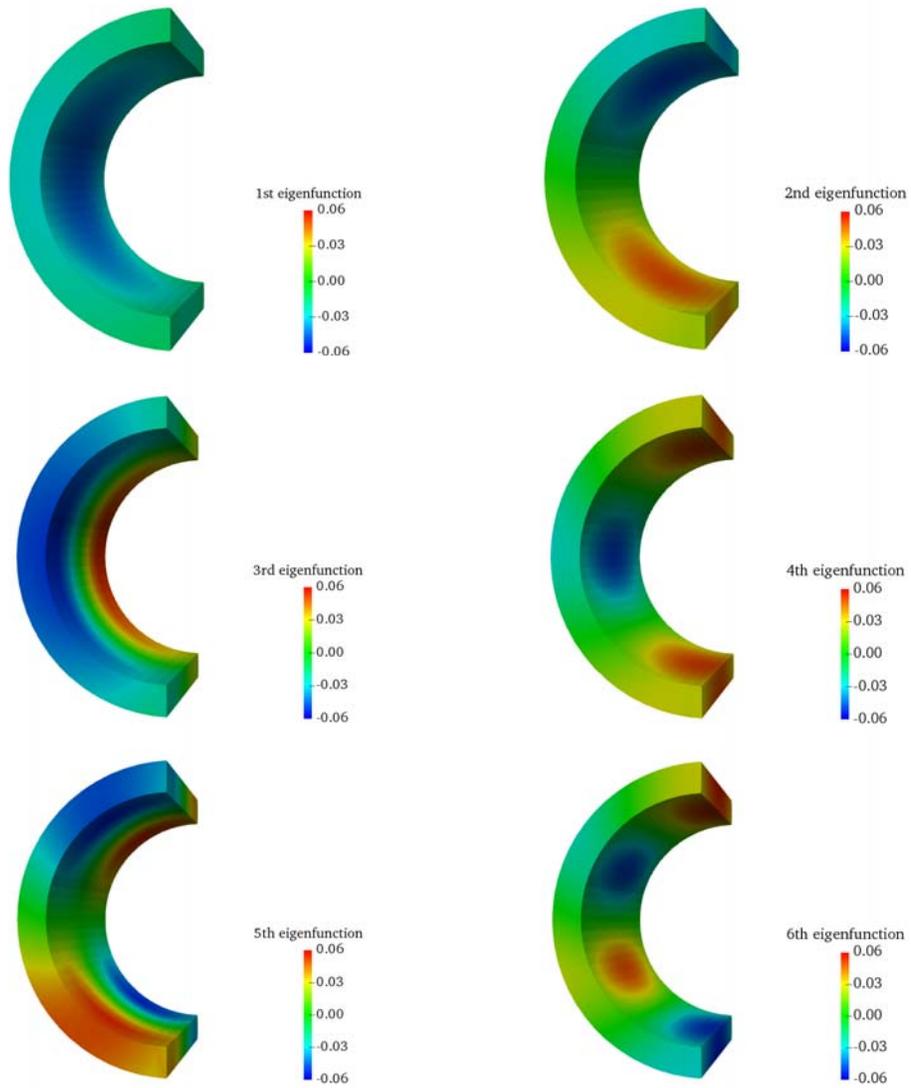

Figure 12: Contour plots of first six IGA-derived eigenfunctions using the fine mesh.

Second, in all three examples presented, the covariance functions and their correlation length parameters are chosen arbitrarily. This is merely for illustration. The method proposed is applicable for any square-integrable covariance function. However, in practice, the covariance function and its correlation length parameter need to be determined from the data collected for the random field problem under consideration. Such an effort is not trivial, as, in many cases, the actual data are noisy. In which case, methods of Bayesian inference and/or statistical estimation techniques will be required to characterize a random field from real-world applications. Therefore, future efforts on how to collect such data and how to estimate the correlation length parameter, including making an optimal choice for the covariance function, should be pursued.

Third, design of complex structures or systems requires estimation of the probabilistic characteristics of an output response variable of interest and the underlying risk of failure. The IGA analysis presented in this work provides a means to parameterize a random field by a finite number of input random variables, which is an important first step. Given a mathematical model of a complex system, a natural extension of the work entails characterizing the discrepancy between model-based simulations and physical reality in terms of the statistical moments, probability law, and other relevant properties of the output variable. For practical applications, encountering hundreds of input random variables is not uncommon, where the output function, defined algorithmically via expensive numerical calculations, is all too often expensive to evaluate. Most existing surrogate methods available today begin to break down when the input-output mapping is highly nonlinear and the input uncertainty is large. More importantly, many high-dimensional problems are all but impossible to solve using existing methods. The root deterrence to practical computability is often related to the high dimension of the multivariate integration or interpolation problem, known as the curse of dimensionality. Therefore, new or improved computational methods capable of exploiting low effective dimensions of multivariate functions are highly desirable.

## 7. Conclusion

A Galerkin isogeometric method was developed for solving an integral eigenvalue problem, resulting in an effective discretization of random fields by means of the well-known Karhunen-Loève expansion. The method employs a Galerkin discretization, which projects the eigensolutions onto a finite-dimensional subspace of a Hilbert space. Using B-splines and NURBS functions as the basis of the subspace, a concomitant matrix eigenvalue problem is formulated, where the system matrices are constructed by domain integrations. Finally, the eigensolutions are obtained using standard methods. Although there exist similar Galerkin methods, such as the finite-element and mesh-free methods, the NURBS-based isogeometric method offers a few computational advantages. First, many physical or computational domains, such as freeform and sculptured surfaces and conic sections, are exactly represented by NURBS. In consequence, potential numerical errors originating from imprecise geometry, accepted in the finite-element and mesh-free methods, are avoided. Second, as NURBS functions have higher-order continuity, the eigensolutions derived from isogeometric analysis are usually globally smoother than those derived from finite-element analysis. The smoothness can be controlled by judiciously selecting or adjusting the polynomial order of the underlying B-splines as well as the multiplicity of knots. Therefore, the introduction of the isogeometric method for random field discretization is not only novel, but it also presents an attractive alternative to existing methods. More importantly, using NURBS for random field discretization enhances the isogeometric paradigm. In consequence, one can envision developing a seamless uncertainty quantification pipeline, where geometric modeling, stress analysis, and stochastic simulation are all consolidated using the same building blocks of NURBS. Numerical results, obtained for three random field discretization problems in all three dimensions, indicate that the isogeometric method developed provides accurate and convergent eigensolutions.

## Appendix A. IGA Details of Numerical Examples

Tables A.1 and A.2 list the control points and weights for the coarsest mesh in Examples 1 and 2, respectively. Table A.3 describes the knot vectors and polynomial orders, including basic mesh properties, for both the coarse and fine meshes in Example 3.

Table A.1: Control points and weights for the coarsest one-element IGA mesh in Example 1.

| Linear elements ($p=1$) | | Quadratic elements ($p=2$) | | Cubic elements ($p=3$) | |
|---|---|---|---|---|---|
| $\mathbf{C_i}$ | $w_i$ | $\mathbf{C_i}$ | $w_i$ | $\mathbf{C_i}$ | $w_i$ |
| (0,0) | 1 | (0,0) | 1 | (0,0) | 1 |
| (1,0) | 1 | (1/2,0) | 1 | (1/3,0) | 1 |
| | | (1,0) | 1 | (2/3,0) | 1 |
| | | | | (1,0) | 1 |

Table A.2: Control points and weights for the coarsest one-element IGA mesh in Example 2.

| $\mathbf{C_i}$ | $w_i$ |
|---|---|
| (0.6,0) | 1 |
| (1,0) | 1 |
| (0.6,0.6) | $\frac{1}{\sqrt{2}}$ |
| (1,1) | $\frac{1}{\sqrt{2}}$ |
| (0,0.6) | 1 |
| (0,1) | 1 |

Table A.3: Knot vectors and polynomial orders for the coarse and fine IGA meshes in Example 3.[a]

| Coarse mesh | Fine mesh |
|---|---|
| $\Xi_1 = (0,0,0,1,1,1)$ | $\Xi_1 = (0,0,0,1,1,1)$ |
| $\Xi_2 = \left(0,0,0,\frac{1}{16},\frac{1}{8},\frac{3}{16},\frac{1}{4},\frac{5}{16},\frac{3}{8},\frac{7}{16},\frac{1}{2},\right.$ $\left.\frac{1}{2},\frac{9}{16},\frac{5}{8},\frac{11}{16},\frac{3}{4},\frac{13}{16},\frac{7}{8},\frac{15}{16},1,1,1\right)$ | $\Xi_2 = \left(0,0,0,\frac{1}{32},\frac{1}{16},\frac{3}{32},\frac{1}{8},\frac{5}{32},\frac{3}{16},\frac{7}{32},\frac{1}{4},\frac{9}{32},\frac{5}{16},\right.$ $\frac{11}{32},\frac{3}{8},\frac{13}{32},\frac{7}{16},\frac{15}{32},\frac{1}{2},\frac{1}{2},\frac{17}{32},\frac{9}{16},\frac{19}{32},\frac{10}{16},\frac{21}{32},\frac{11}{16},$ $\left.\frac{23}{32},\frac{3}{4},\frac{25}{32},\frac{13}{16},\frac{27}{32},\frac{7}{8},\frac{29}{32},\frac{15}{16},\frac{31}{32},1,1,1\right)$ |
| $\Xi_3 = \left(0,0,0,\frac{1}{8},\frac{1}{4},\frac{3}{8},\frac{1}{2},\frac{5}{8},\frac{3}{4},\frac{7}{8},1,1,1\right)$ | $\Xi_3 = \left(0,0,0,\frac{1}{8},\frac{1}{4},\frac{3}{8},\frac{1}{2},\frac{5}{8},\frac{3}{4},\frac{7}{8},1,1,1\right)$ |
| $p_1 = p_2 = p_3 = 2$ | $p_1 = p_2 = p_3 = 2$ |
| No. of elements, $n_e = 128$ | No. of elements, $n_e = 256$ |
| No. of control points, $n_c = 570$ | No. of control points, $n_c = 1050$ |

[a] 1 = radial direction; 2 = circumferential direction; 3 = length direction.

# References

[1] W. Betz, I. Papaioannou, D. Straub, Numerical methods for the discretization of random fields by means of the karhunen-loève expansion, Computer Methods in Applied Mechanics and Engineering 271 (2014) 109–129.
[2] K. Karhunen, Über lineare methoden in der wahrscheinlichkeitsrechnung, Ann. Acad. Sci. Fenn. Ser. A. I. 37 (1947) 3–79.
[3] M. Loève, Fonctions aléatoires de second ordre, in: Processus Stochastic et Mouvement Brownien, Gauthier Villars: Paris.
[4] M. Loève, Probability Theory, Vol II., Springer: Berlin, Heidelberg, New York, 1977.
[5] O. P. Le Maître, O. M. Knio, Spectral Methods for Uncertainty Quantification, Springer: Dordrecht Heidelberg London New York, 2010.
[6] E. I. Fredholm, Sur une classe d'equations fonctionnelles, Acta Mathematica 27 (1903) 365–390.
[7] R. Ghanem, P. D. Spanos, Stochastic finite elements: a spectral approach, World Publishing Corp., 1991.
[8] T. J. R. Hughes, The Finite Element Method: Linear Static and Dynamic Finite Element Analysis, Dover Publications, 2000.
[9] T. Belytschko, Y. Y. Lu, , L. Gu, Element-free galerkin methods, International Journal for Numerical Methods in Engineering 37 (1994) 229–256.
[10] S. Rahman, H. Xu, A meshless method for computational stochastic mechanics, International Journal of Computational Methods in Engineering Science and Mechanics 64 (2005) 41–58.
[11] S. Rahman, Meshfree methods in computational stochastic mechanics, in: Reliability-Based Civil Engineering (edited by A. Haldar), World Scientific Publishing Company: Singapore, 2005.
[12] J. A. Cottrell, T. J. R. Hughes, Y. Bazilevs, Isogeometric Analysis: Toward Integration of CAD and FEA, John Wiley & Sons, 2009.
[13] T. J. R. Hughes, J. A. Cottrell, , Y. Bazilevs, Isogeometric analysis: Cad, finite elements, nurbs, exact geometry and mesh refinement, Computer Methods in Applied Mechanics and Engineering 194 (2005) 4135–4195.
[14] L. A. Piegl, W. Tiller, The NURBS Book, Second Edition, Springer-Verlag: Berlin, 1997.
[15] M. Cox, The numerical evaluation of b-splines (1971).
[16] C. De Boor, On calculation with b-splines, Journal of Approximation Theory 6 (1972) 50–62.
[17] G. E. Farin, NURBS Curves and Surfaces: From Projective Geometry to Practical Use, A. K. Peters, Ltd.: Natick, MA, 1999.
[18] Y. Bazilevs, L. Beirao de Veiga, J. Cottrell, T. J. R. Hughes, G. Sangalli, Isogeometric analysis: approximation, stability and error estimates for h-refined meshes, Mathematical Models and Methods in Applied Sciences 16 (2006) 1031–1090.
[19] N. Dunford, J. T. Schwartz, Linear Operators, Spectral Theory, Self Adjoint Operators in Hilbert Space, Part 2, Wiley-Interscience, 1988.
[20] K. Atkinson, The numerical solution of integral equations of the second kind, Cambridge University Press: Cambridge, UK, 1997.
[21] J. Mercer, Functions of positive and negative type and their connection with the theory of integral equations, Philosophical Transactions of the Royal Society A 209 (1909) 415–446.
[22] C. Schwab, R. A. Todor, Karhunenloève approximation of random fields by generalized fast multipole methods, Journal of Computational Physics 217 (2006) 100–122.
[23] M. Grigoriu, Applied non-Gaussian Processes, Prentice-Hall: Englewood Cliffs, NJ, 1988.
[24] F. Auricchio, F. Calabro, T. J. R. Hughes, A. Reali, G. Sangalli, A simple algorithm for obtaining nearly optimal quadrature rules for nurbs-based isogeometric analysis, Computer Methods in Applied Mechanics and Engineering 249252 (2012) 15–27.
[25] M. Eiermann, O. G. Ernst, E. Ullmann, Computational aspects of the stochastic finite element method, Computing and Visualization in Science 10 (2007) 3–15.
[26] B. N. Khoromskij, A. Litvinenko, H. G. Matthies, Application of hierarchical matrices for computing the karhunenloève expansion, Computing 84 (2009) 49–67.
[27] S. Börm, L. Grasedyck, W. Hackbusch, Introduction to hierarchical matrices with applications, Engineering Analysis with Boundary Elements 27 (2003) 405–422.
[28] R. B. Lehoucq, D. C. Sorensen, C. Yang, ARPACK Users Guide: Solution of Large Scale Eigenvalue Problems with Implicitly Restarted Arnoldi Methods, SIAM: Philadelphia, PA, 1998.
[29] K. Wu, H. Simon, Thick-restart lanczos method for large symmetric eigenvalue problems, SIAM Journal on Matrix Analysis and Applications 22 (2000) 602–616.
[30] L. Greengard, V. Rokhlin, A new version of the fast multipole method for the laplace equation in three dimensions, Acta Numerica 6 (1997) 229–269.
[31] MATLAB, Version 2016a, The MathWorks, Inc.: Natick, MA, 2016.